\documentclass[twoside,11pt]{article}

\usepackage{dsfont}
\usepackage{amsfonts}
\usepackage{amssymb}
\usepackage{amsmath}
\usepackage{stmaryrd}
\usepackage{exscale,cmmib57}
\usepackage{authblk}
\usepackage{enumerate}
\usepackage{empheq}
\usepackage[figuresright]{rotating}
\usepackage{caption}
\usepackage{subcaption}
\usepackage{multirow}
\usepackage{fancyhdr}
\usepackage{mdframed}
\usepackage{appendix}
\usepackage[table]{xcolor}
\usepackage{booktabs}
\usepackage{graphicx}
\usepackage{float}

\usepackage[colorlinks = true,
linkcolor = red,
urlcolor  = blue,
citecolor = blue,
anchorcolor = blue]{hyperref}

\newtheorem{thm}{\bfseries Theorem}[section]
\newtheorem{defi}{\bfseries Definition}[section]
\newtheorem{lem}{\bfseries Lemma}[section]
\newtheorem{prop}{\bfseries Proposition}[section]
\newtheorem{cor}{\bfseries Corollary}[section]
\newtheorem{rem}{\it Remark}[section]

\newcommand{\R}{\mathbb{R}}
\newcommand{\N}{\mathbb{N}}

\newcommand{\Int}{\mathop{\rm int}}

\newcommand{\Min}{\mathop{\rm Min}}

\newcommand{\argmin}{\mathop{\rm argmin}}

\title{\Large\bf Generalized Reduced Jacobian Method}

\author[1]{M. El Maghri\thanks{Corresponding author. Email: elmaghri@yahoo.com}}
\author[2]{Y. Elboulqe\thanks{Email: elboulqeyoussef@gmail.com}}
\affil[1]{\small\em Department of Mathematics and Computer, Faculty of Sciences A{\"\i}n Chock, 
Hassan II University, Casablanca, BP. 5366, Morocco}
\affil[2]{\small\em Department of Mathematics, Faculty of Sciences Semlalia, 
	Cadi Ayyad University, Marrakech, BP. 2390, Morocco}
\date{\today}

\begin{document}

\maketitle

\noindent
{\small

\noindent
{\bf Abstract} \quad 
In a recent work, we presented the reduced Jacobian method (RJM) as an extension of Wolfe's 
reduced gradient method to multicriteria (multiobjective) optimization problems dealing with 
linear constraints. This approach reveals that using a reduction technique of the Jacobian 
matrix of the objective avoids scalarization. In the present work, we intend to generalize 
RJM to handle nonlinear constraints too. In fact, we propose a generalized reduced Jacobian 
(GRJ) method that extends Abadie-Carpentier's approach for single-objective programs. To 
this end, we adopt a global reduction strategy based on the fundamental theorem of implicit 
functions. In this perspective, only a reduced descent direction common to all the criteria 
is computed by solving a simple convex program. After establishing an Armijo-type line search 
condition that ensures feasibility, the resulting algorithm is shown to be globally convergent,  
under mild assumptions, to a Pareto critical (KKT-stationary) point. Finally, experimental 
results are presented, including comparisons with other deterministic and evolutionary 
approaches. \\[.3cm]
{\bf Keywords} \quad  Multicriteria optimization. Pareto optima. Nonlinear programming. Pareto 
KKT-stationarity. Reduced Jacobian method. Nonlinear constraints.\\[.3cm]

\medskip

\section{Introduction}

In many areas, like economics, medicine, design, transportation and so on, we are faced with 
multicriteria (or multiobjective) optimization problems (MOP), where several functions must 
be optimized at the same time. Because of the conflicts which may arise between these functions, 
almost never a single point will minimize or maximize them at once. So the Pareto optimality 
concept have to be considered. A point is said to be Pareto optimum or efficient solution, 
if it cannot be improved with respect to a criterion without degrading at least one of the others. 
Taking account the principle, Pareto solutions are incomparable and no single point can represent 
all the others. The problem then consists of determining the set of all Pareto solutions or their 
objective values called the Pareto front or more generally the set of non-dominated solutions. 

\smallskip

In recent years, it has been observed that this task has become more tractable with the advent 
of a new generation of algorithms grounded in the principle of multiobjective search descent 
directions, which aim to identify a direction that  decreases all the criteria simultaneously 
(e.g. \cite{coc,dru,el2,el4,fl1,fl2, gar,mig,mor,muk,per,qu}). In most cases, the resulting 
algorithms are merely direct extensions of well-known nonlinear descent methods, which tackle 
the given problem without resorting to intermediate transformations or introducing artificial 
parameters that may be sensitive to the original problem. From a theoretical perspective, 
convergence toward Pareto-KKT stationarity is well established, and numerically, these approaches 
have also proven very promising, both in terms of convergence to the true Pareto front and the 
diversity of non-dominated approximate solutions. For instance, in \cite{el2}, we proposed a 
direct extension of the Wolfe reduced gradient method (RGM) \cite{wol} (see also \cite{el1} for 
a new RGM variant) to the so-called reduced Jacobian method (RJM) for solving linearly constrained 
MOPs. We also refer the reader to \cite{el4}, where we introduced an RJM variant that possesses 
the full convergence property.

\smallskip

In the present study, we build upon our previous work by focusing on the so-called generalized 
reduced gradient (GRG) method, originally introduced by Abadie and Carpentier \cite{aba}, which 
extends Wolfe’s reduction strategy to problems with nonlinear constraints. Accordingly, we propose 
a novel approach, which we term the generalized reduced Jacobian (GRJ) method, to effectively 
handle nonlinear constraints in multicriteria optimization. The GRG method is widely recognized 
as one of the fundamental and practical approaches to nonlinear mathematical programming. Its 
convergence properties were rigorously established in subsequent works by Smeers \cite{sme}, 
Mokhtar-Kharroubi \cite{mok}, and Schittkowski \cite{sch}. As far as we are concerned, our 
investigations cover both the theoretical and algorithmic aspects of this method. Abadie-Carpentier’s 
reduction strategy, based on the well-known implicit function theorem, is adopted in this work, 
thereby enabling a complete characterization of Pareto-KKT stationarity through multiobjective 
reduced descent directions. It is subsequently proven that such directions guarantee simultaneous 
descent of all criteria and ensure the existence of Armijo-like steplengths, while preserving 
feasibility throughout the process. From a computational point of view, a simple direction-finding 
subproblem is introduced to iteratively compute the reduced directions. Under mild assumptions, the 
resulting GRJ algorithm is proven to be globally convergent to a Pareto KKT-stationary point in the 
sense of accumulation points. Finally, numerical experiments are reported, illustrating the effectiveness 
of the GRJ method compared with other methods, including a Zoutendijk-like method \cite{mor}, an  
SQP-type approach \cite{fl3}, and the well-known evolutionary method NSGA-II \cite{deb3}.

\medskip

\section{Efficiency and multiobjective descent}

We are dealing with the following nonlinearly constrained multicriteria optimization problem: \\[.1cm]
$$
\text{(MOP)} \qquad \Min\  F(x) \quad \text{subject to} \quad G(x)=0, \;\;  a\leq x\leq b,
$$\\[-.1cm]
where $\big[F,G\big] : x\in\mathbb{R}^n\mapsto \big[F(x),G(x)\big] = \left[\big(f_1(x),\ldots,f_r(x)\big),\big(g_1(x),\ldots,g_m(x)\big)\right] 
\in \mathbb{R}^r \times \mathbb{R}^m$ 
are two vector-valued functions continuously differentiable on the feasible set denoted by $S$, 
and, $a=(a_1,\ldots,a_n),\, b=(b_1,\ldots,b_n) \in \mathbb{R}^n$. For $y=(y_1,\ldots,y_p)\in \R^p$ 
and $y'=(y'_1,\ldots,y'_p)\in \R^p$, $y \leq y'$ (resp. $y < y'$), iff $y_i \leq y'_i$ (resp. 
$y_i < y'_i$) ($\forall i$), and, $y\lneq y'$ means $y\leq y'$ and $y\neq y'$. The usual norm 
and the usual inner product in $\R^p$ will be respectively denoted by $\|y\|=\sqrt{y\cdot y}\,$ 
and $\,y\cdot y'=y^T y'$, where $y^T$ stands for the transpose of $y$.  

\bigskip

The (MOP) consists in minimizing the vector function $F$ over $S$ in the following Pareto senses:

\begin{defi}\label{def1}\em
A point $x^*\in S$ is said to be
\begin{enumerate}[(a)]
\item {\em weakly efficient} for (MOP), if $\ \nexists x\in S, \ F(x)<F(x^*)$;
\item {\em efficient} ({\em Pareto optimum}) for (MOP), if $\ \nexists x\in S, \ F(x)\lneq F(x^*)$;
\item {\em properly efficient} (in Henig's sense) for (MOP), if there exists $C\subsetneq \R^r$ a 
pointed\footnote{A pointed cone $C$ means that $C\cap -C=\left\{0\right\}$.} convex cone such that its 
topological interior $\Int C\supset\R^r_+\backslash\left\{0\right\}$, and, 
$\nexists x\in S, \ F(x)\lneq_C F(x^*)$.\footnote{The relation $y \lneq_C y'$ means that 
$y'-y \in C\backslash\left\{0\right\}$.}
\end{enumerate}
\end{defi}

The sets of weakly efficient points, efficient points and properly efficient points will be denoted 
respectively by $E_w$, $E_e$ and $E_p$. To unify the notations, we briefly denote by $E_{\sigma}$ the 
set of all $\sigma$-efficient solutions depending on the choice of $\sigma\in\left\{w,e,p\right\}$. 
Recall the following relationships:
\begin{equation}\label{rde}
E_p \subseteq E_e \subseteq E_w.
\end{equation}

The concept of {\em local $\sigma$-efficiency} is defined by replacing $S$ by $S\cap{\cal N}(x^*)$ 
in Definition \ref{def1}, where ${\cal N}(x^*)$ is some neighbourhood of $x^*$. The set of all locally 
$\sigma$-efficient points will be denoted by $E_{\sigma}^{loc}$. The image set $F(E_{\sigma})$ (resp. 
$F(E_{\sigma}^{loc})$) is called {\em $\sigma$-Pareto front} (resp. {\em local $\sigma$-Pareto front}). 

\medskip

Taking into account the inclusions (\ref{rde}), the determination of the weakly efficient set will include 
all the other efficient solutions.

\medskip

Let us denote by $JF(x)=\big(\partial f_j (x)/\partial x_i\big)_{j,i}$ the Jacobian matrix of $F$ at $x$; 
the $j$th line of $JF(x)$ being of course the gradient $\nabla f_j(x)$ of the $j$th objective $f_j$ at $x$.

\begin{defi}\label{def2}\em
Let $S$ be a convex set. The vector mapping $F : S \subset \mathbb{R}^n \to \mathbb{R}^r$ is said to be
\begin{enumerate}[(a)]
\item {\em convex} on $S$, if \\[-.15cm]
$$
\forall x, x'\in S,\; \forall \alpha\in \left[0,1\right],\quad F(\alpha x+(1-\alpha)x') \leq \alpha F(x)+(1-\alpha)F(x');
$$		
\item {\em strictly convex} on $S$, if \\[-.1cm]
$$
\forall x, x'\in S,\; x \neq x',\; \forall \alpha\in \left]0,1\right[,\quad F(\alpha x+(1-\alpha)x') < \alpha F(x)+(1-\alpha)F(x');
$$	
\item {\em pseudoconvex} on $S$, if it is pseudoconvex at any $x\in S$, i.e., \\[-.1cm]
$$
\forall x'\in S,\quad F(x')<F(x) \Longrightarrow JF(x)(x'-x) < 0;
$$		
\item {\em strictly pseudoconvex} on $S$, if it is strictly pseudoconvex at any $x\in S$, i.e., \\[-.1cm]
$$
\forall x'\in S,\; x'\neq x,\quad F(x')\leq F(x) \Longrightarrow JF(x)(x'-x) < 0.
$$
\end{enumerate}
\end{defi}

\begin{rem}\label{rem2}\em
It is obvious that $F$ is componentwise convex, iff $F$ is convex. But the converse does 
not hold for pseudoconvexity (see \cite[Theorem 9.2.3 and Remark 5.3]{goh}). The same should 
hold true for the strict concepts. As in the scalar case, pseudoconvexity generalizes convexity. 
\end{rem}

A multiobjective search direction is defined as follows.

\begin{defi}\label{def3}\em
A vector $d\in \R^n$ is said to be
\begin{enumerate}[(a)]
\item a {\em feasible direction} at $x\in S$, if 
$ \ \exists t_f>0, \ \forall t\in \left]0,t_f\right], \ x+td \in S$; \\[-.4cm]
\item a {\em tangent direction} to $S$ at $x\in S$, if 
$\ \exists (d^k)_k \subset \R^n$, 
$\exists (t_k)_k \subset \left]0, +\infty \right[\,$ such that $\lim\limits_{k\rightarrow +\infty} d^k = d$, 
$\lim\limits_{k\rightarrow +\infty} t_k = 0\,$ and $\,x+t_k d^k\in S\,$ ($\forall k$); \\[-.3cm]
\item a {\em descent direction} of $F$ at $x\in S$, if $ \ \exists t_d>0, \ \forall t\in \left]0,t_d\right], \ F(x+td) < F(x)$.
\end{enumerate}
\end{defi}

The sets of feasible directions, tangent directions and descent directions are cones which  
will be denoted respectively by ${\cal A}_S(x)$, ${\cal T}_S(x)$ and ${\cal D}_S(x)$. 

\bigskip

A sufficient condition for $d\in {\cal D}_S(x)$ is (see \cite{fl1})
\begin{equation}\label{eq2}
JF(x)d<0.
\end{equation}

\smallskip

A geometric necessary condition for weak efficiency is obtained using the tangent directions:

\smallskip

\begin{thm} \label{cng}
If $x^*\in E_w^{loc}$, then $\,\nexists d \in {\cal T}_S(x^*)$, $JF(x^*)d<0$.
\end{thm}

\smallskip
\noindent
{\em Proof}\quad Let $d \in {\cal T}_S(x^*)$. We can suppose that $d\neq 0$, otherwise the 
result is trivial. So, by the very definition, 
$\exists (d^k)_k \longrightarrow d, \ \exists (t_k)_k \searrow 0^+$ such that 
\begin{equation}\label{eq4bis}
x^k:=x^*+t_kd^k\in S\quad (\forall k).
\end{equation}
By hypothesis, $\exists\, {\cal N}(x^*)$ a neighbourhood of $x^*$ such that
\begin{equation}\label{eq4}
\nexists x\in {\cal N}(x^*) \cap S, \ F(x)< F(x^*). 
\end{equation}
From (\ref{eq4bis}), we see that $x^k\in{\cal N}(x^*) \cap S$ for $k$ large enough, and 
from (\ref{eq4}), we deduce that
\smallskip 
\begin{equation}\label{eq5}
\exists K\in \N, \ \forall k\geq K, \ \exists j_k\in\left\{1,\ldots,r\right\}, 
\ f_{j_k}(x^k)\geq f_{j_k}(x^*).
\end{equation}
\smallskip 
Without loss of generality, we can assume that $j_k=j_0$ (a constant index) for $k$ large 
enough. On the other hand, by differentiability, we have that 
$$
f_{j_0}\left(x^k\right)=f_{j_0}\left(x^*+t_k d^k\right) = 
f_{j_0}\left(x^*\right)+t_k\nabla f_{j_0}\left(x^*\right)^T d_k+o\left(t_kd^k\right).
$$
It follows, from (\ref{eq5}), that for $k$ large enough, 
$$
t_k\nabla f_{j_0}\left(x^*\right)^T d_k+o\left(t_kd^k\right)\geq 0.
$$
Passing to the limit, when $k\nearrow +\infty$, after dividing by $t_k\|d^k\|$, we obtain
$$
\nabla f_{j_0}\left(x^*\right)^T d\geq 0,
$$
which shows the result of the theorem. \hfill$\hspace{1cm}\Box$

\medskip

As is well known, the tangent cone ${\cal T}_S(x)$  is introduced to handle nonlinear equality 
constraints, since the cone ${\cal A}_S(x)$ may be empty in such cases (see, e.g., the case 
where $S$ is a circle). On the other hand, the tangent cone ${\cal T}_S(x)$  is always closed 
and coincides with the closure of ${\cal A}_S(x)$ when $S$ is convex. This helps to explain the 
suitability of using ${\cal T}_S(x)$ instead of ${\cal A}_S(x)$. Furthermore, although the tangent 
cone plays a theoretically crucial role for a general feasible set $S$, its explicit determination 
often remains challenging in practice. However, when $S$ is explicitly described by equality and/or 
inequality constraints, the tangent cone can always be expressed in terms of these constraints via 
the so-called {\em linearized cone}. For our (MOP), the latter is formulated at a given $x\in S$ 
as follows: 

$$
{\cal C}_S(x)=\left\{d\in \mathbb{R}^n : \ JG(x)d= 0, \ \forall i\in I_a(x),\ d_i\geq 0, 
\ \forall i\in I_b(x),\ d_i\leq 0 \right\},\\[.1cm]
$$

\smallskip 
\noindent
where $I_a(x)=\left\{i \in \left\{1,\ldots,n\right\} :  x_i=a_i\right\}$ and 
$I_b(x)=\left\{i \in \left\{1,\ldots,n\right\} :  x_i=b_i\right\}$ are the index sets of 
active variables, and, $JG(x)$ is the Jacobian matrix of $G$ at $x$. In general, one only 
has that ${\cal T}_S(x) \subseteq {\cal C}_S(x)$ (see, e.g., \cite[Chapter 5]{baz}). When 
equality holds, we say that {\em Abadie's constraint qualification} $({\cal ACQ})$ is 
satisfied at $x$:

\bigskip 

$({\cal ACQ}) \hskip2cm {\cal T}_S(x) = {\cal C}_S(x).$\\
 
\smallskip 

With the assumption $({\cal ACQ})$, we can state the following explicit conditions of weak 
efficiency that will be useful in the sequel.

\begin{thm}\label{thm1}
Assume that $({\cal ACQ})$ holds at $x^*\in S$. If $x^*\in E_w^{loc}$, then the following 
two equivalent conditions are satisfied:
\begin{enumerate}[(i)]
\item There is no tangent descent direction to $S$ at $x^*\in S$ satisfying (\ref{eq2}), i.e.,
\begin{center}
$\nexists d \in \mathbb{R}^n : \quad 
\left\{\begin{aligned}
&JF(x^*)d<0,\\
&JG(x^*)d= 0, \\ 
&d_i\geq 0, \ \forall i\in I_a(x^*),\\
&d_i\leq 0, \ \forall i\in I_b(x^*).
\end{aligned}\right.$
\end{center} 
\item The point $x^*\in S$ is Pareto KKT-stationary, i.e.,
\begin{center}
$\exists \left(\lambda,u,v,w\right)\in \mathbb{R}^r_+\backslash\left\{0\right\} \times \mathbb{R}^m\times\mathbb{R}^{n}_+\times\mathbb{R}^{n}_+ : \quad 
\left\{\begin{aligned}
&JF\left(x^*\right)^T\lambda+JG\left(x^*\right)^Tu-v+w=0,\\
&v\cdot (x^*-a)=0,\\
&w\cdot (b-x^*)=0.
\end{aligned}\right.$
\end{center}
\end{enumerate}
Conversely, if condition {\em (i)} or {\em (ii)} is satisfied at  $x^*\in S$, $F$ is 
pseudoconvex (resp. strictly pseudoconvex) at $x^*$ and the kernel $\text{Ker}\,JG(x^*)$ 
contains ${\cal N}(0)\cap (S-x^*)$ for some neighbourhood ${\cal N}(0)$ of $0$, then 
$x^*\in E_w^{loc}$ (resp. $x^*\in E_e^{loc}$). 
\end{thm}

\smallskip

\noindent
{\em Proof}\quad \textit{Necessity:} condition (i) is the result of Theorem \ref{cng} 
expressed under $({\cal ACQ})$. Equivalence between (i) and (ii) comes directly from 
Motzkin's alternative theorem (see, e.g., \cite{man}).\\
\textit{Sufficiency:} Since it is assumed that ${\cal N}(0)\cap (S-x^*)\subseteq\text{Ker}\,JG(x^*)$, 
then there exists ${\cal N}(x^*)$ a neighbourhood of $x^*$ such that $JG(x^*)(x-x^*)=0$ for 
all $x\in {\cal N}(x^*)\cap S$. Now, by way of contraposition, suppose that $x^*\not\in E^{loc}_w$ 
(resp. $x^*\not\in E^{loc}_e$). Then, there exists $x^0\in S\cap{\cal N}(x^*)$ such that 
$F(x^0)<F(x^*)$ (resp. $F(x^0)\lneq F(x^*)$). By putting $d=x^0-x^*$, pseudoconvexity (resp. 
strict pseudoconvexity) of $F$ at $x^*$ implies that $JF(x^*)d=JF(x^*)\left(x^0-x^*\right)<0$. 
On the other hand, we also have that $JG(x^*)d=JG(x^*)(x^0-x^*)=0$ because 
$x^0\in {\cal N}(x^*)\cap S\subseteq\text{Ker}\,JG(x^*)$. Moreover,
$$
d_i=\left\{\begin{aligned}
&x^0_i-a_i, \ \text{if} \ i\in I_a(x^*),\\ 
&x^0_i-b_i, \ \text{if} \ i\in I_b(x^*),
\end{aligned}\right.
$$\\
satisfying $d_i\geq 0$, if $i\in I_a(x^*)$ and $d_i\leq 0$, if $i\in I_b(x^*)$, which under 
the constraint qualification condition $({\cal ACQ})$, shows the existence of a tangent descent 
direction of $F$ at $x^*$. \hfill$\hspace{1cm}\Box$

\begin{rem}\label{rem3}\em
If, in Theorem \ref{thm1}, we replace ${\cal N}(0)$ by $\R^n$, then it is easy to see by the 
same proof that the local efficiency becomes global, in particular, when $G$ is affine.
\end{rem}
 
\smallskip

\section{GRJ strategy}

Let $A(x):=JG(x)\in \R^{m \times n}$ be the Jacobian matrix of $G$ at $x\in S$ such that 
$A(x)$ is of full rank $m<n$. Then, there exists a subset $B\subset \left\{1,\ldots,n\right\}$, 
called {\em basis} of $A(x)$, such that the submatrix $A_B(x)$ is invertible, where (rearranging 
the columns if necessary) $A(x)=[A_B(x) \ A_N(x)]$ with $A_B(x)=JG_B(x)$, $A_N(x)=JG_N(x)$ and 
obviously $N=\left\{1,\ldots,n\right\} \backslash B$. Rearranging also the components of $x$ we 
can write $x=(x_B,x_N)$, where $x_B$ (resp.  $x_N$) is the well-known vector of {\em basic} 
(resp. {\em nonbasic}) variables. Since $G$ is assumed to be continuously differentiable and 
$A_B(x)$ is invertible, then by the implicit function theorem, there exist $V\subset \mathbb{R}^{m}$ 
a neighbourhood of $x_B$, $W\subset \mathbb{R}^{n-m}$ a neighbourhood of $x_N$ and a unique function 
$\psi : \ W \longrightarrow V$ which is continuously differentiable such that 
\begin{eqnarray}\label{eq6}
\forall x'_N \in W, \qquad 
\begin{cases}
G\left(\psi(x'_N),x'_N\right)=0,\\[0.2cm]
\displaystyle{\frac{\partial \psi}{\partial x_N}}(x'_N)=-A_{B}^{-1}(x')A_N(x').
\end{cases}
\end{eqnarray}\\ 
Hence, the feasible set $S$ may be expressed only in terms of nonbasic variables:
$$
x\in S \Longleftrightarrow  a\leq x=\left(\psi(x_N),x_N\right)\leq b.
$$

A point $x\in S$ is said to be {\em nondegenerate}, if there exists $B$ a basis such that 
$a_B<x_B<b_B$. In this case, $B$ is also said to be a nondegenerate basis for $x$; otherwise 
$x$ is {\em degenerate}. The feasible set $S$ is in turn called nondegenerate, if any 
$x\in S$ is nondegenerate. 

\medskip

Given a vector $d\in \R^n$ that we partition in the same way to $d=(d_B,d_N)$, then it is 
clear that
$$
JG(x)d=0 \ \Longleftrightarrow \ d_B=-A_{B}^{-1}(x)A_N(x)d_N. 
$$
\smallskip 
So by also partitioning the Jacobian matrix $JF(x)=[JF_B(x) \ JF_N(x)]$, for any $d$ such 
that $d_B=-A_{B}^{-1}(x)A_N(x)d_N$, it follows that
$$ 
JF(x)d=U_N(x)d_N,
$$
where $U_N(x)$ is an $r\times(n-m)$ matrix, which will be called {\em generalized reduced 
Jacobian matrix} of $F$ at $x$, explicitly given by
\begin{equation}\label{rjm}
U_N(x):=JF_N(x)-JF_B(x)A_{B}^{-1}(x)A_N(x).
\end{equation}
Observe that the $j$th row of $U_N(x)$ is nothing more than the well-known reduced gradient 
of the criterion $f_j$ at $x$ which is given by 
$$
u^j_N(x) = 
\frac{\partial f_j}{\partial x_N}(x)-\left( A_{B}^{-1}(x)A_N(x)\right) ^T\frac{\partial f_j}{\partial x_B}(x).
$$ 

Now, we define a {\em multiobjective reduced descent direction} of $F$ at $x$ as being a 
nonbasic vector $d_N\in \R^{n-m}$ satisfying:
\begin{equation}\label{eq7}
U_N(x)d_N<0, \quad \forall i\in I_a(x)\cap N,\ d_i\geq 0, \quad \text{and} \quad  
\forall i\in I_b(x)\cap N,\ d_i\leq 0.
\end{equation}

\medskip

The geometric stationarity stated in Theorem \ref{thm1}(i) then can be expressed only in 
terms of nonbasic directions:

\begin{cor} \label{cor1} 
Under $({\cal ACQ})$, a nondegenerate point $x\in S$ is Pareto KKT-stationary for (MOP), iff 
there is no multiobjective reduced descent direction of $F$ at $x$.	
\end{cor}

The next result shows that given any multiobjective reduced descent direction at a given 
nondegenerate point, feasibility along this direction with an Armijo-like steplength is 
guaranteed.

\begin{prop}\label{fac}
Let $x\in S$ and $\psi$ be the associated implicit function given by (\ref{eq6}), and, let 
$d_N$ be a multiobjective reduced descent direction of $F$ at $x$. Let $t\in \mathbb{R}_+$ 
and define
\begin{equation}\label{tN}
t_N:=\min\left(\left\lbrace\frac{a_i-x_i}{d_i} \ : \ d_i<0, \ i\in N\right\rbrace \cup 
\left\lbrace\frac{b_i-x_i}{d_i} \ : \ d_i>0, \ i\in N \right\rbrace\right). 
\end{equation}\\[-.3cm]
Then, 
\begin{enumerate}[(i)]
\item $t_N>0$. 

\item $t\leq t_N \Longleftrightarrow a_N\leq x_N + td_N \leq b_N$.
\end{enumerate}

\smallskip 
\noindent
If furthermore $a_B < x_B <b_B$ (i.e., $x$ is nondegenerate with $B=\{1,\ldots,n\}\setminus N$), 
then

\begin{enumerate}[(i)]
\item[(iii)] $\exists t_f\in \left]0, t_N\right]$, $\forall t\in \left]0, t_f\right]$,  
\begin{equation}\label{aa}
\big(\psi(x_N + td_N),x_N + td_N\big) \in S\quad \text{and} \quad 
a_B < \psi(x_N + td_N) < b_B.
\end{equation}
\item[(iv)] $\forall \beta\in \left]0,1 \right[$ (fixed), $\exists t_a\in \left]0, t_f\right]$, 
$\forall t\in \left]0, t_a\right]$,
\begin{equation}\label{ac}
F\big(\psi(x_N + td_N),x_N + td_N\big) <F(x)+\beta t U_N(x)d_N.
\end{equation}
\end{enumerate}
\end{prop}

\smallskip

\noindent
{\em Proof}\quad
(i) It is clear that $t_N$ exists (since $d_N\neq 0$). The positivity of $t_N$ follows by 
the very definition of $d_N$. Indeed, if $t_N=(a_i-x_i)/d_i$, then from (\ref{eq7}), it 
follows that $a_i-x_i<0$, and consequently, $t_N>0$; otherwise, $t_N=(b_i-x_i)/d_i$, which 
in accordance with (\ref{eq7}), entails that $b_i-x_i>0$, thereby ensuring $t_N>0$.\\[0.2cm]
(ii) Let us show the first implication. Let $i\in N$. The assertion is obvious if $d_i=0$. 
Now, if $d_i>0$, then $a_i\leq x_i+td_i$ and $(b_i-x_i)/d_i\geq t_N\geq t$, so that 
$b_i\geq x_i+td_i$. The case $d_i<0$ is similar. Conversely, since $d_N\neq 0$, then if 
$d_i<0$ ($i\in N$), the hypothesis $a_i\leq x_i+td_i$ implies that $t\leq(a_i-x_i)/d_i$; 
and if $d_i>0$, then the hypothesis $x_i+td_i\leq b_i$ implies that $t\leq(b_i-x_i)/d_i$, which 
shows that $t\leq t_N$.\\[0.2cm]
(iii) By definition of $\psi : \ W \longrightarrow V$ (see (\ref{eq6})) and the continuity 
of $t\longmapsto x_N + td_N$ at $0$, there exists $t_1>0$ such that $x_N + td_N\in W$ for 
all $t\in\left[0,t_1\right]$. It follows by (\ref{eq6}) that 
$$
\forall t\in\left[0,t_1\right], \quad G\left(\psi(x_N + td_N),x_N + td_N\right)=0.
$$
On the other hand, by assertion (ii), we have that 
$$
\forall t\in\left[0,t_{N}\right], \quad a_N\leq x_N + td_N\leq b_N.
$$
So, it remains to show that the box-constraint also holds for the basic variables. Indeed, 
by nondegeneracy assumption and continuity, it is immediate that 
$$
\exists t_2>0, \quad \forall t\in\left[0,t_2\right], \quad a_B < \psi(x_N + td_N) < b_B.
$$
Thus, by putting $t_f=\min\left\{t_1,t_2,t_{N}\right\}$, we obtain the desired result.\\[0.2cm]
(iv) To prove the Armijo-like inequality (\ref{ac}), let us consider the reduced function 
$\widetilde{F}(\cdot):=F\left( \psi(\cdot),\cdot\right)$ which is continuously differentiable 
on $W\subset \mathbb{R}^{n-m}$. By using the derivative of $\psi$ given by (\ref{eq6}), we 
immediately deduce the Jacobian of $\widetilde{F}$:
$$
J\widetilde{F}(x_N)=JF_N(x)-JF_B(x)A_B(x)^{-1}A_N(x)=U_N(x).
$$  
By virtue of (iii), $x_N+td_N\in W$  for all $t\in \left]0, t_f\right]$. So, by differentiability, 
it follows that 
$$
\begin{array}{lcl}
\widetilde{F}(x_N+td_N)&=&\widetilde{F}(x_N)+tU_N(x)d_N+o(t)\\[0.1cm]
                       &=&\widetilde{F}(x_N)+t\beta U_N(x)d_N+ t\big[(1-\beta) U_N(x)d_N +\frac{1}{t}o(t)\big].
\end{array}
$$
Since $(1-\beta)U_N(x)d_N<0$ and $\lim\limits_{t\rightarrow 0^+}\dfrac{1}{t}o(t)=0$, then \\
$$
\exists t_a\in \left]0, t_f\right],\ \forall t\in \left]0, t_a\right],\;\;(1-\beta) U_N(x)d_N +\frac{1}{t}o(t)<0,
$$
thus implying that
$$
\widetilde{F}(x_N+td_N)<\widetilde{F}(x_N)+t\beta U_N(x)d_N.
$$\\[-.1cm]
As we have $\widetilde{F}(x_N+td_N)=F\big(\psi(x_N+td_N),x_N+td_N\big)$ and $\psi(x_N)=x_B$, 
then the result follows straightforwardly.
\hfill$\hspace{1cm}\Box$

\medskip 

After having shown that it is possible to obtain feasible steplengths along reduced directions, 
it remains for completing our GRJ scheme to show how to determine such reduced descent directions.

\medskip

\section{GRJ search direction}

Let us consider first a positive functional $\phi : \R \rightarrow \R_+$ satisfying 
$\phi(t)=0$ iff $t=0$; e.g., $\phi_p(t)=\displaystyle{\frac{\left|t\right|^p}{p}}$ 
for some value of $p\in \left]0,1\right]$, or, $\phi_0(t)=\mathds{1}_{\R^*}$ the 
characteristic function of $\R^*$ defined by $\mathds{1}_{\R^*}(t)=1$ if $t\neq0$; 
$\mathds{1}_{\R^*}(0)=0$. Let $\left\lfloor a \right\rfloor_+ := \max(0, a)$ and 
$\left\lfloor a \right\rfloor_- := \max(0, -a)$ denote, respectively, the positive 
and negative parts of the scalar $a$. Then, for $x\in S$, we introduce the following 
finding-direction subproblem:
$$
({\cal P}_x) \qquad \underset{\lambda\in \Lambda}{\Min} \ \ f(\lambda, x):= 
\displaystyle{\frac{1}{2}}\sum\limits_{i\in N}\left(\phi\left(b_i-x_i\right)\Big\lfloor 
\left(U_N(x)^T\lambda\right)_i\Big\rfloor^2_- + \phi\left(x_i-a_i\right)\Big\lfloor 
\left(U_N(x)^T\lambda\right)_i \Big\rfloor^2_+ \right)  
$$\\
where $\Lambda=\{ (\lambda_1,\ldots,\lambda_r)\in \R^r_+ \ : \ \sum_{j=1}^{r} \lambda_j=1\}$, 
$x_i$ (resp. $a_i$ and $b_i$) is the $ith$ component of $x_N$ (resp. $a_N$ and $b_N$) and 
$\left(U_N(x)^T\lambda\right)_i$ is the $ith$ component of $U_N(x)^T\lambda$. 

\medskip

It is obvious that $({\cal P}_x)$ is a convex continuously differentiable problem and always 
has an optimal solution in the compact set $\Lambda$. Moreover, $({\cal P}_x)$ has simplicial 
constraints and it is easy to verify that its objective function has a computable gradient 
explicitly given by   
$$ 
\nabla_{\lambda}f(\lambda,x)=-U_N(x)\delta_N(\lambda,x), 
$$
where 
\begin{equation} \label{delta}
\delta_N(\lambda,x)=\left(\phi\left(b_i-x_i\right)\Big\lfloor \left(U_N(x)^T\lambda\right)_i\Big\rfloor_- 
- \phi\left(x_i-a_i\right)\Big\lfloor \left(U_N(x)^T\lambda\right)_i \Big\rfloor_+\right)_{i\in N}.
\end{equation} 

\medskip

The descent GRJ scheme will consist to take $d_N:=\delta_N(\lambda(x),x)$, where 
$\lambda(x)\in \argmin({\cal P}_x)$. Then, according to (\ref{delta}), $\forall i\in N$, 
\begin{equation}\label{rdd}
d_i=\left\{
\begin{aligned}
&-\phi(x_i-a_i)\left(U_N(x)^T\lambda(x)\right)_i, &\;& \text{if}\quad \left(U_N(x)^T\lambda(x)\right)_i>0, \\[0.2cm]
&-\phi(b_i-x_i)\left(U_N(x)^T\lambda(x)\right)_i, &\;&\text{else}. 
\end{aligned}
\right.
\end{equation}

\medskip

The lemma below characterizes the optimal solutions of $({\cal P}_x)$, which will be 
crucial not only for proving the descent properties of the vector $d_N$, but also for 
the convergence analysis of the resulting algorithm. Its proof is omitted, as it is 
similar to the one of \cite[Lemma 3.2]{el2}.

\begin{lem} \label{lac} $\lambda(x)\in \Lambda$ is an optimal solution of $({\cal P}_x)$, 
iff
$$
\delta_N(\lambda(x),x)\cdot\left( U_N(x)^T\lambda\right)\leq -2f\left(\lambda(x),x\right) 
\qquad (\forall \lambda \in \Lambda),
$$
or equivalently,
$$
\delta_N(\lambda(x),x)\cdot u^j_N(x)\leq -2f\left(\lambda(x),x\right) \qquad (j=1,\ldots,r).
$$
\end{lem}

\medskip 

The descent properties are stated as follows:

\begin{prop}\label{mdp} \ 

\begin{enumerate}[(i)]
\item $d_N=0 \Longleftrightarrow f(\lambda(x),x)=0$.

\item $d_N\neq 0 \Longrightarrow d_N$ is a multiobjective reduced descent direction 
of $F$ at $x$. If moreover $({\cal ACQ})$ holds and $x$ is nondegenerate, then $x$ is 
not Pareto KKT-stationary for (MOP). 
		
\item $d_N=0 \Longrightarrow x$ is Pareto KKT-stationary for (MOP). 

\item The functional $x \mapsto f(\lambda(x),x)$ is continuous on $S$ provided $\phi$ 
is continuous.
\end{enumerate}
\end{prop}

\smallskip

\noindent
{\em Proof}\quad (i) To show the direct implication, observe that as $\phi\geq 0$, 
then $d_N=0$ implies that, $\forall i\in N$,
$$
\Big\lfloor \phi\left(b_i-x_i\right) \left(U_N(x)^T\lambda(x)\right)_i\Big\rfloor_- = 
\Big\lfloor \phi\left(x_i-a_i\right) \left(U_N(x)^T\lambda(x)\right)_i\Big\rfloor_+.
$$
If we suppose that the last two terms are not equal to zero, we would have 
$$
\phi\left(b_i-x_i\right) \left(U_N(x)^T\lambda(x)\right)_i <0 \quad \text{and} \quad \phi\left(x_i-a_i\right) \left(U_N(x)^T\lambda(x)\right)_i>0,
$$
which would mean that $\phi\left(b_i-x_i\right)$ and $\phi\left(x_i-a_i\right)$ have 
not the same sign contradicting $\phi\geq 0$. 
Thus,  $\forall i\in N$,
$$
\Big\lfloor \phi\left(b_i-x_i\right) \left(U_N(x)^T\lambda(x)\right)_i\Big\rfloor_- = 
\Big\lfloor \phi\left(x_i-a_i\right) \left(U_N(x)^T\lambda(x)\right)_i\Big\rfloor_+=0, 
$$
and this shows that $f(\lambda(x),x)=0$. The reverse implication follows straightforwardly 
from the positivity of $\phi$ and  (\ref{rdd}).\\[0.2cm]
(ii) If $d_N\neq0$, then by assertion (i), $f(\lambda(x),x)>0$. It follows, by Lemma \ref{lac}, 
that
$$
d_N\cdot u^j_N(x)\leq -2f(\lambda(x),x)<0 \qquad (j=1,\ldots,r).
$$
Hence, $U_N(x)d_N<0$. On the other hand, for $i\in I_a(x)\cap N$, we have that
$$
d_i=\phi(b_i-x_i)\Big\lfloor\left(U_N(x)^T\lambda(x)\right)_i\Big\rfloor_- \geq0,
$$ 
and, for $i\in I_b(x)\cap N$, we have that 
$$
d_i=-\phi(a_i-x_i)\Big\lfloor\left(U_N(x)^T\lambda(x)\right)_i\Big\rfloor_+ \leq0. 
$$
This shows, according to the definition (\ref{eq7}), that $d_N$ is well a multiobjective 
reduced descent direction. The second part of this assertion follows straightforwardly 
from Corollary \ref{cor1}. \\[0.2cm]
(iii) Suppose that $d_N=0 $. If we put 
$$
v^*_N:=\left(\Big\lfloor\left(U_N(x)^T\lambda(x)\right)_i\Big\rfloor_+\right)_{i\in N} 
\quad \text{and}\qquad 
w^*_N:=\left( \Big\lfloor\left(U_N(x)^T\lambda(x)\right)_i\Big\rfloor_-\right)_{i\in N},
$$ 
then, by (\ref{rdd}) and the property that $\phi(t)=0$ iff $t=0$, we obtain
$$
(v_N^*, w_N^*)\geq 0, \quad v_N^*\cdot (x_N-a_N)=0 \quad \text{and} \quad w_N^*\cdot (b_N-x_N)=0.
$$
Clearly, for $\lambda^*:=\lambda(x)$, $u^*:=-(A^{-1}_B(x))^T JF_B(x)^T\lambda^*$, $v^*:=(0,v^*_N)$ 
and $w^*:=(0,w^*_N)$, it holds that 
$\left(\lambda^*,u^*,v^*,w^*\right) \in \mathbb{R}^r_+\backslash\left\{0\right\}\times\mathbb{R}^m\times\mathbb{R}^{n}_+\times\mathbb{R}^{n}_+$,  $JF(x)^T\lambda^*+JG(x)^Tu^*-v^*+w^*=0$ and $v^*\cdot (x-a)=w^*\cdot(b-x)=0$, which express 
the Pareto KKT-stationarity for (MOP).\\[0.2cm]
(iv) Since $f(\lambda(x), x)$ is the optimal value of $({\cal P}_x)$, its continuity 
with respect to $x$ follows from the compactness of the feasible set $\Lambda$ and 
the continuity of the function $(\lambda, x)\mapsto f(\lambda,x)$. \hfill$\hspace{1cm}\Box$

\medskip

\section{GRJ algorithm and its convergence}

\subsection{The algorithm}

The generalized reduced Jacobian algorithm we propose to solve the multicriteria problem (MOP) is 
described below:

\begin{center}
\hrule
\vspace{0.3cm}
GRJ pseudocode with Armijo line search
\vspace{0.3cm}
\hrule
\end{center}
{\footnotesize
{\bf Step 0: Initialization.}\\[0.15cm]
- Select $x\in S$, fix the Armijo's constant $\beta\in\left]0,1\right[$ and choose the functional $\phi$.\\[0.15cm]
{\bf Step 1: Nondegenerate basis selection.}\\[0.15cm]
- Identify a basis  $B$ such that $a_B<x_B<b_B$. Set $N=\left\{1,\ldots,n\right\}\backslash B$.\\[0.15cm]
{\bf Step 2: Generalized reduced Jacobian.}\\[0.15cm]
- Compute $U_N(x)$ according to (\ref{rjm}).\\[0.15cm]
{\bf Step 3: Reduced direction.}\\[0.15cm]
- Solve $({\cal P}_x)$ to determine $d_N$ according to (\ref{rdd}).\\[0.15cm]
{\bf Step 4: Stopping criterion.}\\[0.15cm]
- If $\min({\cal P}_x)=0$, then STOP: $x$ is Pareto KKT-stationary for (MOP).\\[0.15cm]
{\bf Step 5: Feasible Armijo line search.}\\[0.15cm]
- Compute $t_N$ according to (\ref{tN}).\\[0.15cm]
- Determine, according to (\ref{aa})-(\ref{ac}), a steplength $t\in \left]0,t_N\right]$ that satisfies: 
$$
{\small t=\underset{p\in \mathbb{N}}{\max}\left\{\frac{1}{2^p} \quad
: \quad \begin{aligned}
&\Big(\psi\big(x_N+\frac{1}{2^p}d_N\big),x_N+\frac{1}{2^p}d_N\Big)\in S, \ \text{and}\\[0.2cm]
&F\Big(\psi\big(x_N+\frac{1}{2^p}d_N\big),x_N+\frac{1}{2^p}d_N\Big)<F(x)+\beta \frac{1}{2^p}U_N(x)d_N
\end{aligned}
\right\rbrace },
$$
where $\psi$ is the implicit function defined in (\ref{eq6}).\\[0.15cm]
{\bf Step 6: Updated point.}\\[0.15cm]
- Set $x:=\left(\psi\left(x_N+td_N\right),x_N+td_N\right)$.\\[0.15cm]
{\bf Step 7: Degeneracy test.}\\[0.15cm]
- If $B$ is a degenerate basis for $x$, go to Step 1;
otherwise go to Step 2.\\
\hrule}

\subsection{Convergence analysis}

If the algorithm terminates after a finite number of iterations, then according to 
Proposition \ref{mdp}, it terminates at a Pareto KKT-stationary point. So we will 
assume that the sequence produced by the algorithm is infinite. We also assume that 
the bases are chosen in such a way that the following hypothesis 
is satisfied:

\medskip 

$$
({\cal H})\qquad \left\{\begin{array}{l}
x^k\underset{k\in K \subseteq\mathbb{N}}{\longrightarrow} x^* \ \text{non Pareto KKT-stationary},\\[.5cm]
\exists  B \ \text{basis}, \ \forall k\in K, \ a_B < x^k_B < b_B.
\end{array}
\right\} \Longrightarrow  a_B < x^*_B < b_B.
$$

\medskip 

Known as the \emph{basis property}, this hypothesis proves to be essential for the 
global convergence of reduced gradient methods (e.g., \cite{el2, el4, mok, sme}
\footnote{As pointed out in \cite{el4}, the hypothesis $({\cal H})$ was inadvertently 
omitted in the proof of \cite[Theorem 4.1]{el2}.}). An example provided 
in \cite{mok} shows that, without this hypothesis, the algorithm may fail to converge 
even if the constraints are linear. However, there are basis selection procedures that 
ensure the basis property for both linear and nonlinear constraint cases, such as those 
mentioned in \cite{el4}.

\begin{thm}\label{cgt} Assume that $S$ is nondegenerate, $\phi$ is chosen to be continuous 
and derivable at $0$. Let $(x^k)_{k\in \N}$ be the sequence produced by the GRJ algorithm. 
Then,
\begin{enumerate}[(i)]
\item $x^k\in S$ and $F(x^{k+1})<F(x^k)$ for all $k\in \N$;
\item under the hypothesis (${\cal H}$), any accumulation point $x^*$ of $(x^k)_{k}$ is a 
Pareto KKT-stationary point for (MOP), and, $\displaystyle{\lim_{k\rightarrow +\infty}F(x^k)=F(x^*)}$.
\end{enumerate}
\end{thm}

\smallskip

\noindent
{\em Proof}\quad (i) This assertion follows directly from Proposition \ref{fac}. \\
\noindent
(ii) We need to prove the following lemma.
\begin{lem}\label{lemc}
If, in addition to the hypotheses of Theorem \ref{cgt}, a subsequence $(x^k)_{k\in K}$ 
produced by the GRJ algorithm converges to $x^*$ with constant bases $B_k=B$ for all 
$k\in K$, and $\left(\lambda(x^k)\right)_{k\in K}$ converges too, then 
\begin{enumerate}[(a)]
\item $d_N^k :=\delta_N(\lambda(x^k),x^k)\underset{k\in K}{\longrightarrow} d_N^* :=\delta_N(\lambda(x^*),x^*)$, 
where $N=\left\{1,\ldots,n\right\}\backslash B$;
\item $\bar{t}:=\underset{k\in K}{\inf} \ t_N^k>0$, where
$$
\ t_N^k:=\min\left(\left\lbrace\frac{a_i-x^k_i}{d_i^k} \ : \ d_i^k<0, 
\ i\in N\right\rbrace \cup \left\lbrace\frac{b_i-x^k_i}{d_i^k} \ : \ d_i^k>0, 
\ i\in N \right\rbrace\right).
$$
\end{enumerate}
\end{lem}

\medskip

\noindent
{\em Proof} \quad (a) Let $\lambda^*=\lim\limits_{k\in K \atop k\to +\infty}\lambda(x^k)$, 
then $\lambda^*\in \Lambda$ since $\left(\lambda(x^k)\right)_{k\in K}\subset \Lambda$ which 
is closed. According to the hypothesis and Lemma \ref{lac}, for all $k\in K$, we have 
$$
\delta_N(\lambda(x^k),x^k)\cdot u^j_N(x^k)\leq -2f(\lambda(x^k),x^k) \qquad (j=1,\ldots,r),
$$
Taking into account that the mappings $(\lambda,x)\longmapsto f(\lambda,x)$ and 
$(\lambda,x)\longmapsto \delta_N(\lambda,x)$ are continuous, by passing onto the 
limit, when $k\nearrow +\infty$ in the previous inequality, we obtain that
$$
\delta_N(\lambda^*,x^*)\cdot u^j_N(x^*)\leq -2f(\lambda^*,x^*)  \qquad (j=1,\ldots,r),
$$
Since by closeness of $S$, $x^*$ remains feasible as $(x^k)_k$, this shows applying 
once again Lemma \ref{lac} that $\lambda^*\in \argmin({\cal P}_{x^*})$, i.e., 
$\lambda^*=\lambda(x^*)$, which proves the assertion. \\[0.2cm]
(b) To prove this assertion, we shall proceed by a contradiction way by assuming 
that $\bar{t}=0$ (because $\bar{t}\geq 0$). Then, it would exist $K_0\subseteq K$ 
an infinite subset such that $t_N^k\underset{k\in K_0}{\longrightarrow}0$. Since 
the index set $N$ is finite, we can assume (without loss of generality) that one of 
the two following cases would happen: 
$$
\exists i_0\in N, \quad \forall k\in K_0, \quad t_N^k=\frac{a_{i_0}-x^k_{i_0}}{d_{i_0}^k}, 
\quad d_{i_0}^k<0,
$$ 
or
$$
\exists i_0\in N, \quad \forall k\in K_0, \quad t_N^k=\frac{b_{i_0}-x^k_{i_0}}{d_{i_0}^k}, 
\quad d_{i_0}^k>0.
$$
In the first case, following (\ref{rdd}), we would have 
\begin{equation}\label{eq14}
\forall k\in K_0, \quad t_N^k=\frac{1}{\left(U_N(x^k)^T\lambda(x^k)\right)_{i_0}}\left[\frac{x^k_{i_0}-a_{i_0}}{\phi(x^k_{i_0}-a_{i_0})}\right].
\end{equation}

\smallskip 
\noindent 
Now, by assertion (a), $d_{i_0}^k := - \phi(x^k_{i_0}-a_{i_0})\left(U_N(x^k)^T\lambda(x^k)\right)_{i_0} 
\underset{k\in K_0}{\longrightarrow} d_{i_0}^*$, and, by the contradiction hypothesis,  
$t_N^k \underset{k\in K_0}{\longrightarrow} 0$. So we would have that the sequence 
$\left(x^k_{i_0} - a_{i_0}\right)_{k \in K_0}$ also converges to 0; i.e., 
$x^k_{i_0} \underset{k\in K_0}{\longrightarrow} a_{i_0}$. But this would imply that 
$$
\phi(x^k_{i_0}-a_{i_0})\underset{k\in K_0}{\longrightarrow}\phi(x^*_{i_0}-a_{i_0})=\phi(0)=0,
$$ 
and, by letting $k\nearrow +\infty$ in (\ref{eq14}), we would obtain the contradiction  
$$
t_N^k\underset{k\in K_0}{\longrightarrow} \frac{1}{\left(U_N(x^*)^T\lambda^*\right)_{i_0}}\times \frac{1}{\phi'(0)}\neq 0,
$$
The second case is similar. Hence, $\bar{t}>0$. \hfill$\hspace{1cm}\Box$

\bigskip 

Let us go back to the proof of Theorem \ref{cgt}(ii). 

\medskip 

By hypothesis, there exists a subsequence 
$(x^k)_{k\in K}\longrightarrow x^*$ with $K\subseteq \mathbb{N}$ an infinite subset. By 
closeness of $S$, $x^*$ remains feasible as $(x^k)_k$. Now, by continuity of $F$ and the 
assertion (i), it follows that the sequence $\big(F(x^k)\big)_{k\in\N}\longrightarrow F(x^*)$. 
On the other hand, since the index set $B_k$ is in the finite set $\{1,\ldots,n\}$, we can 
assume (without loss of generality) that $B_k=B$ (hence $N_k=N$) for all $k\in K$. Now, 
by Armijo's condition: $\forall k\in \N$,
$$
F(x^{k+1})-F(x^k)<\beta t_kU_{N_k}(x^k)d_{N_k}^k<0.
$$
It follows that
\begin{equation}\label{lim}
\underset{k\rightarrow +\infty}{\underset{k\in K}{\lim}}t_kU_N(x^k)d_N^k = 
\underset{k\rightarrow +\infty}{\underset{k\in K}{\lim}}t_kU_{N_k}(x^k)d_{N_k}^k=0.
\end{equation}
Similarly, as the sequence $\left(\lambda(x^k)\right)_{k}\subset \Lambda$ which is a 
compact set, we can assume (without loss of generality) that 
$\lambda(x^k)\underset{k\in K}{\longrightarrow} \lambda^*\in\Lambda$. Then, by virtue of 
Lemma \ref{lemc}(a),  $d_N^k\underset{k\in K}{\longrightarrow} d_N^*$. If we suppose, by 
a contradiction way, that $x^*$ is not Pareto KKT-stationary, then according to 
Proposition ($\ref{mdp}$)(ii)--(iii), $d_N^*$ would be a multiobjective reduced descent 
direction of $F$ at $x^*\in S$. Hence, according to the definition (\ref{eq7}),  
\begin{equation}\label{ds*}
U_N(x^*)d_N^* < 0.
\end{equation}
We can see from (\ref{lim}) that this would imply $t_k \underset{k\in K}{\longrightarrow} 0$.  
In particular, for any $p\in\mathbb{N}$, $t_k<\frac{1}{2^p}$ for any $k\in K$ sufficiently large. 
Recall that, by the implicit function theorem applied at $x^*$, there exists a unique mapping 
$\psi^*: \ W^*\longrightarrow V^*$, where $W^*$ (resp. $V^*$) is a neighbourhood of $x^*_N$ 
(resp. $x^*_B=\psi^*(x^*_N)$) such that $G\big(\psi^*(x_N),x_N\big)=0$ for all $x_N\in W^*$. 
This applies to each $x^k_N+\frac{1}{2^p}d_N^k$ as $(k,p)\to +\infty$, since 
$x^k_N+\frac{1}{2^p}d_N^k \longrightarrow x^*_N$; hence, for all $p\in \N$ and $k\in K$, both 
sufficiently large,  
\begin{equation} \label{adm1}
G\Big(\psi^*\big(x^k_N+\frac{1}{2^p}d_N^k\big),x^k_N+\frac{1}{2^p}d_N^k\Big)=0.
\end{equation}
On the other hand, since $x^*$ is assumed to be not Pareto KKT-stationary, then using the 
hypothesis (${\cal H}$),  $x^*$ would be nondegenerate, so that $a_B < x^*_B=\psi^*(x^*_N) < b_B$. 
Thus, by the continuity of $\psi^*$, we would also have that, for all $p\in \N$ and $k\in K$, 
both sufficiently large, 
\begin{equation} \label{adm2}
a_B < \psi^*\Big(x^k_N+\frac{1}{2^p}d_N^k\Big) < b_B. 
\end{equation}
Now, using Lemma \ref{lemc}(b), we have that $0<\frac{1}{2^p}\leq \bar{t}$ for all $p$ large 
enough, and this shows, according to Proposition \ref{fac}(ii), that for any $p\in \N$ 
sufficiently large and $k\in K$,  
\begin{equation} \label{adm3}
a_N\leq x^k_N+\frac{1}{2^p}d_N^k \leq b_N. 
\end{equation} 
Subsequently, we would have, following (\ref{adm1})--(\ref{adm3}), that  
for all $p\in \N$ and $k\in K$, both sufficiently large,
$$
\Big(\psi^*\big(x^k_N+\frac{1}{2^p}d_N^k\big),x^k_N+\frac{1}{2^p}d_N^k\Big)\in S.
$$

\smallskip
\noindent 
However, given the previous assertion that $t_k < \frac{1}{2^p}$, this would mean that, 
in Step 5 of the GRJ algorithm, for all $p\in \N$ and $k\in K$, both sufficiently large, 
the Armijo inequality would fail to hold; that is,
\begin{equation}\label{eq17}
F\Big(\psi^*\big(x^k_N+\frac{1}{2^p}d_N^k\big),x^k_N+\frac{1}{2^p}d_N^k\Big)\not<F(x^k)+\beta \frac{1}{2^p}U_N(x^k)d_N^k.
\end{equation}

\medskip
\noindent 
Fixing $p\in \N$ sufficiently large and passing onto the limit when $k\nearrow +\infty$, 
we would obtain that
$$
F\Big(\psi^*\big(x^*_N+\frac{1}{2^p}d_N^*\big),x^*_N+\frac{1}{2^p}d_N^*\Big)\not<F(x^*)+\beta \frac{1}{2^p}U_N(x^*)d_N^*.
$$

\smallskip
\noindent 
Since $p\in \N$ was arbitrary large enough, this would be in contradiction with 
Proposition \ref{fac}(iv), and this proves that $x^*$ is well Pareto KKT-stationary 
for (MOP). \hfill$\hspace{1cm}\Box$

\medskip

\section{Numerical experiments}

\subsection{Test problems and implementation}

In this section, we present computational results obtained through the proposed GRJ 
method. We also provide comparisons with three other multicriteria methods: ZMO 
(a Zoutendijk-like method \cite{mor}), MOSQP (an SQP-type approach \cite{fl3}) and 
the evolutionary NSGA-II method (the nondominated sorting genetic algorithm \cite{deb3}).
All codes were implemented on a machine equipped with 1.90 GHz Intel(R) Core(TM) i5 CPU 
and 16 Go memory. GRJ, ZMO, and MOSQP were coded in SCILAB-2024.1.0, while NSGA-II was 
executed using MATLAB's predefined function ``gamultiobj'' (MATLAB R2024b).

\smallskip

Thirty constrained multicriteria optimization test problems, selected from the literature, 
were considered in this experiment and are summarized in Table \ref{T1}, where `OV', `L' 
and `NL', respectively, indicate the number of original variables, the number of linear 
constraints and the number of nonlinear constraints. We selected the `Disc Brake' and 
`Welded Beam' problems, which are derived from real-world engineering design applications 
\cite{ray}. Detailed descriptions of these two models, along with further analysis of the 
results, are provided at the end of this section. Also, we introduced the problem `EL3' 
below showing the importance of the numerical aspect of the tangent cone ${\cal T}_S(x)$, 
in the case where no feasible direction exists, i.e.,  ${\cal A}_S(x) = \emptyset$. 
\begin{align}\label{EL3}
(\text{EL3}) \qquad&\text{Minimize} \quad \left(x_2^3+\log\left(x_1^2+1\right),\sin\Big(\displaystyle{\frac{x_1}{x_2+2}}\Big)\right),\notag\\
&\text{subject to}  \quad x_1^2+x_2^2=1,\\ 
& \qquad \qquad \quad  \ 0\leq x_i \leq 1, \quad i=1,2.\notag
\end{align}

\medskip

\begin{table}[H]
\centering{\footnotesize
\begin{tabular}{l c c c c c|l c c c c c }
Problem & Ref. & r & OV & L & NL&Problem & Ref. & r & OV & L& NL\\ 
\hline
ABC\_comp&\cite{hwa}&2&2&2&1            &LAP1&\cite{coc}&2&2&1&2\\
BNH&\cite{deb1}&2&2&0&2                 &LAP2&\cite{coc}&2&30&0&1\\
CF11&\cite{zho}&2&3&0&2                 &LIR-CMOP1&\cite{fan}&2&30&0&2\\
CPT1&\cite{deb2}&2&2&0&30               &LIR-CMOP2&\cite{fan}&2&30&0&2\\
Disc Brake&\cite{ray}&2&4&2&3           &LIR-CMOP3&\cite{fan}&2&30&0&3\\
DTLZ0&\cite{deb4}&3&3&0&2               &liswetm&\cite{ley}&2&7&5&0\\
DTLZ9&\cite{deb4}&3&30&0&2              &MOLPg\_001&\cite{ste}&3&8&8&0\\
EL3&(\ref{EL3})&2&2&0&1                 &MOLPg\_002&\cite{ste}&3&12&13&0\\
ex001&\cite{das}&2&5&1&2                &MOLPg\_003&\cite{ste}&3&10&12&0\\
ex002&\cite{wan}&2&5&0&4                &OSY&\cite{deb1}&2&6&4&2\\
ex003&\cite{deb1}&2&2&0&2               &SRN&\cite{deb1}&2&2&1&1\\
GE1&\cite{eic}&2&2&0&1                  &Tamaki&\cite{ray}&3&3&0&1\\
GE4&\cite{eic}&3&3&0&1                  &TLK1&\cite{tha}&2&2&4&0\\
Hanne4 &\cite{col}&2&2&0&1              &TNK&\cite{deb1}&2&2&0&2\\
hs05x&\cite{hoc}&3&5&6&0                &Welded Beam&\cite{ray}&2&4&1&3\\                                           
\hline               
\end{tabular}
\caption{Linearly and nonlinearly constrained multicriteria test problems}\label{T1}}
\end{table}

All these problems were solved starting with the same population of $200$ selected 
individuals. To ensure that all compared solvers generate the same number of 200 
final solutions, and since MATLAB's NSGA-II solver may provide fewer than 200, we 
set ``PopulationSize'' option to 571 for NSGA-II, while all other NSGA-II options 
were left with the default settings. 

\smallskip

In our implementation, the three codes GRJ, ZMO and MOSQP use the value 0.25 as 
Armijo’s constant, and
$$
\min\left({\cal P}_x\right) < 10^{-6}
$$
as the stopping criterion, where ${\cal P}_x$ denotes the direction-finding subproblem 
specific to each of the three solvers. The latter was solved, for GRJ, by coding a 
standard RGM based on Wolfe’s continuous scheme (see, e.g., \cite{el4} for details), 
while for ZMO and MOSQP, it was solved using SCILAB’s predefined function ``qld'' for 
linear quadratic programming problems.

\smallskip

Step 5 of the GRJ algorithm, which involves computing feasible Armijo's steplengths, relies 
on the results presented in Proposition \ref{fac}. This step is carried out using the 
Newton method, following the procedure outlined below: 

\begin{center}
\hrule
\vspace{0.3cm}
Feasible Armijo line search pseudocode
\vspace{0.3cm}
\hrule
\end{center}
{\footnotesize
{\bf Step 0: Initialization.}\\[0.15cm]
- Fix a tolerance $\varepsilon>0$ (e.g., $\varepsilon=10^{-6}$) and 
a maximum number of iterations $L\in\N$ (e.g., $L=200$). \\[0.15cm]
- Set  $t=t_N$, $y^1_B=x_B$ and $l=1$. \\[0.15cm]
{\bf Step 1: Newton rule.}\\[0.15cm]
- Set $y^{l+1}_B=y^l_B-A_{B}^{-1}\Big(y^l_B,x_N+td_N\Big)G\Big(y^l_B,x_N+td_N\Big)$. \\[0.15cm]
{\bf Step 2: Feasible Armijo test.}\\[0.15cm]
- If
$$
\left\|G\Big(y^{l+1}_B,x_N+td_N\Big)\right\|<\varepsilon, \quad a_B\leq y^{l+1}_B\leq b_B, \quad \text{and} \quad F\Big(y^{l+1}_B,x_N+td_N\Big)<F(x)+\beta t U_Nd_N,
$$
then STOP: the current implicit basic vector is set as $\psi\left(x_N+td_N\right) = y^{l+1}_B$.\\[0.15cm]
{\bf Step 3: Update.}\\[0.15cm]
- If $l=L$, set $t=\dfrac{t}{2}$, $y^1_B=x_B$ and $l=1$; otherwise, set $l=l+1$.\\[0.15cm]
- Return to Step 1. \\
\hrule}

\subsection{Performance measures and interpretation}

To thoroughly analyse and compare the results obtained by the solvers, a set of 
performance metrics, as suggested in the literature (see, e.g., \cite{aud}), may 
be employed. These metrics essentially assess two key aspects: the convergence of 
the approximate solutions towards the Pareto front and their diversity. In this work, 
we have chosen three performance measures: {\em Purity} (P) \cite{ban}, {\em Spread} 
($\Delta^*$) \cite{zhou} and {\em Generational distance} (GD) \cite{van}. These 
three measures use, as a proxy for the true Pareto front, the so-called 
{\em reference Pareto front},  defined as: \\[-.2cm]
$$
F_p=\left\lbrace y^*\in\underset{s\in {\cal S}}{\bigcup} F_{p,s} \ : \ \nexists  y\in\underset{s\in {\cal S}}{\bigcup} F_{p,s}, \ y<y^*\right\rbrace, 
$$
where $F_{p,s}$ stands for the approximated Pareto front of problem $p\in {\cal P}$ 
provided by method $s\in {\cal S}$; ${\cal P}$ denotes the set of tested problems and 
${\cal S}$ the set of considered solvers.
\begin{description}
\item[$\mathbf{-}$] \textbf{Purity} (P). This metric consists in measuring the proportion 
of points in $F_{p,s}$ admitted in $F_p$:
$$
P_{p,s}=\frac{\left| F_{p,s}\cap F_p\right|}{\left|F_{p,s}\right|}.
$$
Clearly $P_{p,s}\in [0, 1]$ and the extreme values are significant in the sense that a 
value $P_{p,s}$ close to 1 indicates better performance, while a value equal to 0 implies 
that the solver is unable to generate any point of $F_p$.
\item[$\mathbf{-}$] \textbf{Spread} ($\Delta^*$). This metric measures the diversity 
and the dispersion of the points in $F_{p,s}$ with respect to $F_p$:
$$
\Delta^*_{p,s}=\frac{\sum\limits_{j=1}^{r} d\left(y_j^*, F_{p,s}\right) 
+\sum\limits_{y\in F_p}\big|d\left(y, F_{p,s}\setminus\{y\}\right) - \bar{d}\big|}
{\sum\limits_{j=1}^{r} d\left(y_j^*, F_{p,s}\right) + \left|F_{p}\right|\bar{d}},
$$
where $d(y,A)=\min_{a\in A}\left\| y-a\right\|$ is the Euclidean distance from the vector 
$y$ to the set $A$;  here, $y_j^*$ denotes the minimum value of the criterion $f_j$ in $F_p$, 
and $\bar{d}$ represents the average distance between each solution $y \in F_p$ and the set 
$F_{p,s} \setminus \{y\}$. Note that lower values of $\Delta^*_{p,s}$ reflects more uniform 
distributions of the generated solutions. 
\item[$\mathbf{-}$] \textbf{Generational Distance} (GD). Measuring the convergence, this 
metric represents how far $F_{p,s}$ is from $F_p$:
$$
GD_{p,s}=\frac{1}{\left|F_{p,s}\right|}\sqrt{\sum\limits_{y\in F_{p,s}} d^2\left(y, F_p\right) }.
$$
Obviously, lower values of $GD_{p,s}$ are requested.
\end{description}

Table \ref{T2} lists the numerical results obtained by the four considered methods. 
The column `CPU' indicates the average execution time (in seconds). The best scores 
are in bold. To analyse these results and look for significant differences between 
the compared solvers, the {\em performance profile} is used (see, e.g., \cite{dol}). 
Recall that a profile is a graphical representation of the (cumulative) distribution 
function $\rho$ whose outputs are a solver's scores on the set of test problems against 
a performance metric. More precisely, given the measure $m_{p,s}$ (set here to 1/P, 
$\Delta^*$, GD or CPU) by solver $s\in {\cal S}$ for solving $p\in {\cal P}$, the 
corresponding distribution function is given formally by
$$
\rho_s=\frac{\left|\left\lbrace p\in{\cal P} \ : \ r_{p,s}\leq \alpha \right\rbrace\right|}{\left|{\cal P}\right|},
$$
where $r_{p,s} = m_{p,s}/ \min \left\lbrace m_{p,s} : s\in S\right\rbrace $. The purity 
measure $m_{p,s}$ is set to $1/P$ in order that all the considered metrics exhibit the 
same asymptotic behaviour, in a sense that the smaller the measure, the better the solver. 
Note that at the threshold $\alpha=1$, $\rho_s(\alpha)$ gives us the largest number of 
problems among the best solved by $s$ according to the analysed performance. However, 
a value $\rho_s(\alpha)$ attaining 1 means that all the problems $p\in {\cal P}$ have 
been solved by solver $s$ at the threshold $\alpha$. Thus, the best overall performance 
of a solver is that which reaches the value 1 for the smallest value of $\alpha$. 

\begin{table}[H]
\centering{\footnotesize
\vspace*{-0.9cm}
\begin{tabular}{l@{\hskip0.3em}c@{\hskip0.3em}l@{\hskip0.5em}l@{\hskip0.5em}l@{\hskip0.6em}l@{\hskip0.7em}c@{\hskip0.7em}l@{\hskip0.5em}l@{\hskip0.5em}l@{\hskip0.3em}l@{\hskip-0.2em}}
	\toprule
	\multirow{2}{4em}{Problems}&&\multicolumn{4}{l}{\hskip-0.2cm GRJ}&&\multicolumn{4}{l}{\hskip-0.2cm ZMO}\\
	\cline{3-6}\cline{8-11}
	&&CPU               &P              &$\Delta^\star$    &GD                     && CPU              &P              &$\Delta^\star$    &GD\\
	\toprule
	ABC\_comp    &&\textbf{0.015935} &0.99           &\textbf{0.902345} &0.0000076              &&0.021454          &\textbf{1}     &0.953210          &\textbf{0}\\
	BNH          &&1.436124          &0.915          &0.810000          &0.0000086              &&\textbf{0.004504} &\textbf{1}     &0.995432          &\textbf{0}\\
	CF11         &&0.415047          &0.885          &0.920000          &\textbf{0}             &&0.058855          &0.61           &\textbf{0.902345} &\textbf{0}\\
	CPT1         &&1.712818          &0.97           &\textbf{0.765432} &0.000003               &&0.055902          &0.98           &0.895432          &0.000003\\
	Disc Brake   &&1.242267          &0.69           &\textbf{0.310123} &0.0031797              &&1.817904          &0.035          &0.860987          &0.004630\\
	DTLZ0        &&0.410519          &\textbf{0.805} &0.705678          &0.000615               &&\textbf{0.000827} &0             &0.765432          &0.001699\\			
	DTLZ9        &&\textbf{0.209932} &0.9            &0.254321          &\textbf{0}             &&1.004320          &0.08           &0.989765          &0.000092\\   
	EL3          &&\textbf{0.005883} &\textbf{1}     &\textbf{0.285432} &\textbf{0}             &&0.005889          &\textbf{1}     &0.925678          &\textbf{0}\\
	ex001        &&0.112501          &\textbf{1}     &0.415098          &\textbf{0}             &&0.020554          &0.33           &\textbf{0.350987} &0.276242\\
	ex002        &&0.027517          &0.46           &\textbf{0.785432} &0.0027654              &&\textbf{0.010278} &0              &0.825678          &0.042952\\
	ex003        &&0.044351          &0.535          &0.811023          &0.0003054              &&\textbf{0.038164} &\textbf{0.565} &0.985432          &0.000682\\
	GE1          &&0.018449          &\textbf{0.995} &0.870987          &\textbf{0}             &&\textbf{0.002708} &\textbf{0.995} &0.990987          &\textbf{0}\\
	GE4          &&0.048618          &0.985          &0.632098          &0.0000041              &&0.037197          &\textbf{0.99}  &0.991098          &0.000004\\
	Hanne4       &&0.024853          &0.56           &\textbf{0.707098} &0.0004533              &&0.030180          &0.46           &0.998098          &0.000614\\
	hs05x        &&1.008168          &0.975          &0.980987          &0.0006024              &&\textbf{0.006076} &0.01           &0.970987          &0.019143\\
	LAP1         &&1.822010          &0.755          &\textbf{0.505432} &\textbf{0.000108}      &&0.028570          &\textbf{0.965} &0.640987          &0.000128\\
	LAP2         &&0.087737          &0.87           &\textbf{0.180123} &0.0058377              &&0.082840          &\textbf{1}     &0.995432          &\textbf{0}\\
	LIR-CMOP1    &&0.389043          &\textbf{0.837} &0.840987          &0.000219               &&\textbf{0.035736} &0              &0.965098          &0.003135\\
	LIR-CMOP2    &&\textbf{0.135677} &\textbf{0.9}   &\textbf{0.750987} &0.000038               &&0.135736          &0.02           &0.997453          &0.003135\\
	LIR-CMOP3    &&0.031187          &\textbf{0.995} &\textbf{0.509871} &\textbf{0}             &&\textbf{0.010200} &0.62           &0.740987          &0.000065\\
	liswetm      &&0.154443          &\textbf{0.61}  &\textbf{0.150987} &\textbf{0.000048}      &&\textbf{0.015055} &0              &0.930987          &0.000815\\
	MOLPg\_001   &&\textbf{0.184877} &0.91           &0.955432          &0.000004               &&0.550163          &0.81           &0.980987          &0.000005\\
	MOLPg\_002   &&1.200010          &0.625          &0.740987          &0.000758               &&\textbf{0.015298} &0              &0.720987          &0.040700\\
	MOLPg\_003   &&1.005378          &\textbf{0.99}  &0.700987          &\textbf{0}             &&0.042475          &0              &0.560987          &0.021736\\
	OSY          &&1.539581          &0.47           &0.923444          &0.003442               &&\textbf{0.061214} &0              &\textbf{0.781341} &0.019456\\
	SRN          &&1.937231          &\textbf{0.945} &0.860987          &\textbf{0}             &&\textbf{0.010315} &0.9            &0.940987          &0.006855\\
	Tamaki       &&1.541475          &0.83           &0.700987          &\textbf{0}             &&\textbf{0.002469} &\textbf{1}     &\textbf{0.530987} &\textbf{0}\\
	TLK1         &&\textbf{0.013301} &0.42           &\textbf{0.840987} &0.000101               &&0.033947          &0.39           &0.950097          &0.000080\\
	TNK          &&0.459437          &\textbf{0.395} &0.830987          &0.0003194              &&\textbf{0.013540} &0.275          &\textbf{0.760987} &0.004407\\
	Welded Beam  &&1.137963          &0.91           &\textbf{0.740987} &0.0033562              &&13.20645          &0              &0.960987          &1.226983\\
	\toprule
	\multirow{2}{4em}{Problems} &&\multicolumn{4}{l}{\hskip-0.2cm MOSQP}                        &&\multicolumn{4}{l}{\hskip-0.2cm NSGA-II}\\
	\cline{3-6}\cline{8-11}
	&&CPU               &P              &$\Delta^\star$    &GD                    &&CPU                 &P              &$\Delta^\star$    &GD\\
	\toprule
	ABC\_comp    &&0.020315          &0.99           &0.980987          &0.000017              &&0.857500            &0.41           &0.970000          &0.000906\\
	BNH          &&0.163576          &0.105          &0.992109          &0.010358              &&0.013100            &0.84           &\textbf{0.780000}          &0.001654\\
	CF11         &&\textbf{0.011990} &0.59           &1.000000          &\textbf{0}            &&0.021800            &\textbf{1}  &0.994567          &\textbf{0}\\
	CPT1         &&\textbf{0.017768} &0.74           &0.900000          &0.000025              &&0.191800            &\textbf{1}     &1.000000          &\textbf{0}\\
	Disc Brake   &&0.021562          &0              &0.360000          &0.006729              &&\textbf{0.0084}     &\textbf{0.84}  &0.840012          &\textbf{0.001951}\\
	DTLZ0        &&0.014875          &0.495          &0.960000          &\textbf{0.000528}     &&0.0094              &0.485          &\textbf{0.610000 }&0.003421\\
	DTLZ9        &&0.395620          &\textbf{0.93}  &0.840000          &0.005980              &&0.3825431           &0.08           &\textbf{0.070981} &0.000092\\
	EL3          &&0.012636          &0              &0.659765          &0.010208              &&0.0143              &0.99           &0.591220          &0.000028\\
	ex001        &&\textbf{0.087176} &0              &1.000000          &1.875443              &&4.9669              &0              &0.991197          &\textbf{0}\\
	ex002        &&0.029159          &\textbf{0.55}  &0.967787          &\textbf{0.001515}     &&0.038400            &0.085          &1.000000          &0.021188\\
	ex003        &&0.065472          &0.395          &\textbf{0.793245} &0.000364              &&0.147100            &0       &0.985023          &\textbf{0.000061}\\
	GE1          &&0.031198 &\textbf{0.995} &1.000000          &\textbf{0}            &&0.010200            &0.62    &\textbf{0.760000} &0.000040\\
	GE4          &&\textbf{0.019889} &0.97           &0.992501          &\textbf{0}            &&0.01200             &0.87           &\textbf{0.570000} &0.000458\\
	Hanne4       &&0.022001          &0.355          &0.980000          &\textbf{0.000308}     &&\textbf{0.014400}   &\textbf{0.775} &0.900000          &0.000643\\
	hs05x        &&0.056837          &\textbf{1}     &0.990000          &\textbf{0}            &&3.059300            &0.941          &\textbf{0.905401} &0.010439\\
	LAP1         &&0.023713          &0.255          &0.970350          &0.000203              &&\textbf{0.010900}   &0.78           &0.870000          &0.000176\\
	LAP2         &&0.024061          &0.27           &0.790000          &0.006000              &&\textbf{0.016500}   &0.075          &0.960000          &0.059720\\
	LIR-CMOP1    &&0.059540          &0.8            &\textbf{0.806701} &0.000053              && 0.1824925          &0.22    &1.000000          &\textbf{0.000002}\\
	LIR-CMOP2    &&0.236534          &0.53           &0.950276          &\textbf{0.000022}     &&2.1357355           &0.1            &0.993204          &0.041562\\
	LIR-CMOP3    &&0.323865          &0.76           &0.805467          &0.000376              &&\textbf{0.010200}   &0.62           &0.760000          &0.000065\\
	liswetm      &&0.016368          &0.385          &1.000000          &0.000410              &&1.119700            &0.515          &0.830000          &0.000244\\
	MOLPg\_001   &&0.1888712         &\textbf{0.96}  &\textbf{0.934124} &0.000234              &&0.210634            &0.9            &0.993204          &\textbf{0}\\
	MOLPg\_002   &&0.017967          &\textbf{0.68 } &0.995119          &\textbf{0}            &&1.95200             &0.075          &\textbf{0.703005} &0.075590\\
	MOLPg\_003   &&\textbf{0.032949} &0.79           &0.968578          &0.001239              &&1.693900            &0.19           &\textbf{0.503204} &0.042700\\
	OSY          &&1.582346          &\textbf{0.525} &0.998123          &\textbf{0.001645}     &&1.837200            &0.25           &0.934551          &0.001930\\
	SRN          &&0.036814          &\textbf{0.945} &0.967634          &\textbf{0}            &&2.393800            &0.685          &\textbf{0.773497}          &0.004229\\
	Tamaki       &&0.02418           &0.9            &0.843457          &\textbf{0}            &&0.017700            &0.99           &0.952078          &0.000007\\
	TLK1         &&0.019948          &0.225          &0.996780          &0.000042              &&0.959700            &\textbf{0.985} &0.990840          &\textbf{0.000004}\\
	TNK          &&1.014764          &0.09           &1.000000          &\textbf{0.000245}     &&0.97000             &0.33           &1.000000          &0.003962\\
	Welded Beam  &&\textbf{0.024213} &0.05           &0.999732          &0.178934              &&0.0875              &\textbf{0.97}  &0.950000          &\textbf{0.000621}\\
	\hline
\end{tabular}
\caption{Multiobjective performance measurements}\label{T2}}
\end{table}

\medskip 
 
As illustrated in Fig.~\ref{fig1}, the GRJ method consistently demonstrates 
good performance with respect to the P metric, compared to the other methods.
Indeed, the best value $\rho(\alpha) = 1$ is achieved by GRJ at a relatively 
small threshold $\alpha$, whereas the other methods fail to reach this value. 
This observation is further supported by the results in Tab.~\ref{T2}, where 
the purity values reported for the ZMO, MOSQP and NSGA-II methods are in some 
cases very close or equal to zero. This suggests that these methods fail to 
reach the reference Pareto front for certain problems. 

\smallskip

In terms of the $\Delta^*$ and GD metrics, the four methods appear to be generally 
comparable. Although NSGA-II is widely recognized for its effectiveness in spread, 
the observed profile indicates good overall performance across all solvers, with a 
slight advantage noted for GRJ and NSGA-II. Similar observations can be made regarding 
the GD metric, noting this time a slight advantage of the three methods over NSGA-II, 
which is well known for its limited convergence toward the Pareto front.

\smallskip

Regarding CPU time, as shown in Fig.~\ref{fig1} and Tab.~\ref{T2}, the three 
methods exhibit a slight speed advantage over GRJ. This can be attributed to the 
inherent characteristics of the GRJ algorithm and its operational mechanism. More 
specifically, the algorithm’s iterative process and heuristic strategies, such as 
the basis selection procedure in Step 1 and the line search procedure in Step 5 
(as previously described), can be time-consuming. Nevertheless, these components 
enable efficient exploration of the solution space, which helps to explain its 
performance differences relative to the other methods. The graphical representation 
of the approximated Pareto fronts illustrated in Fig.~\ref{fig2}--\ref{fig3} also 
supports these interpretations. Note also that ZMO generally appears to be faster 
than MOSQP based on the scores given in Tab.~\ref{T2}, except for the case where 
ZMO relatively took a long time to solve the 'Welded Beam' design problem. This 
may explain the influence of MOSQP's CPU performance profile on ZMO observed in 
Fig.~\ref{fig1}. However, as reported in Tab.~\ref{T2}, the execution time remains 
generally very small and reasonable for all the considered solvers. 

\begin{figure}[ht!]
	\centering
	\hspace*{-0.8cm}
	\begin{subfigure}[b]{0.46\textwidth}
		\includegraphics[height=2.4in]{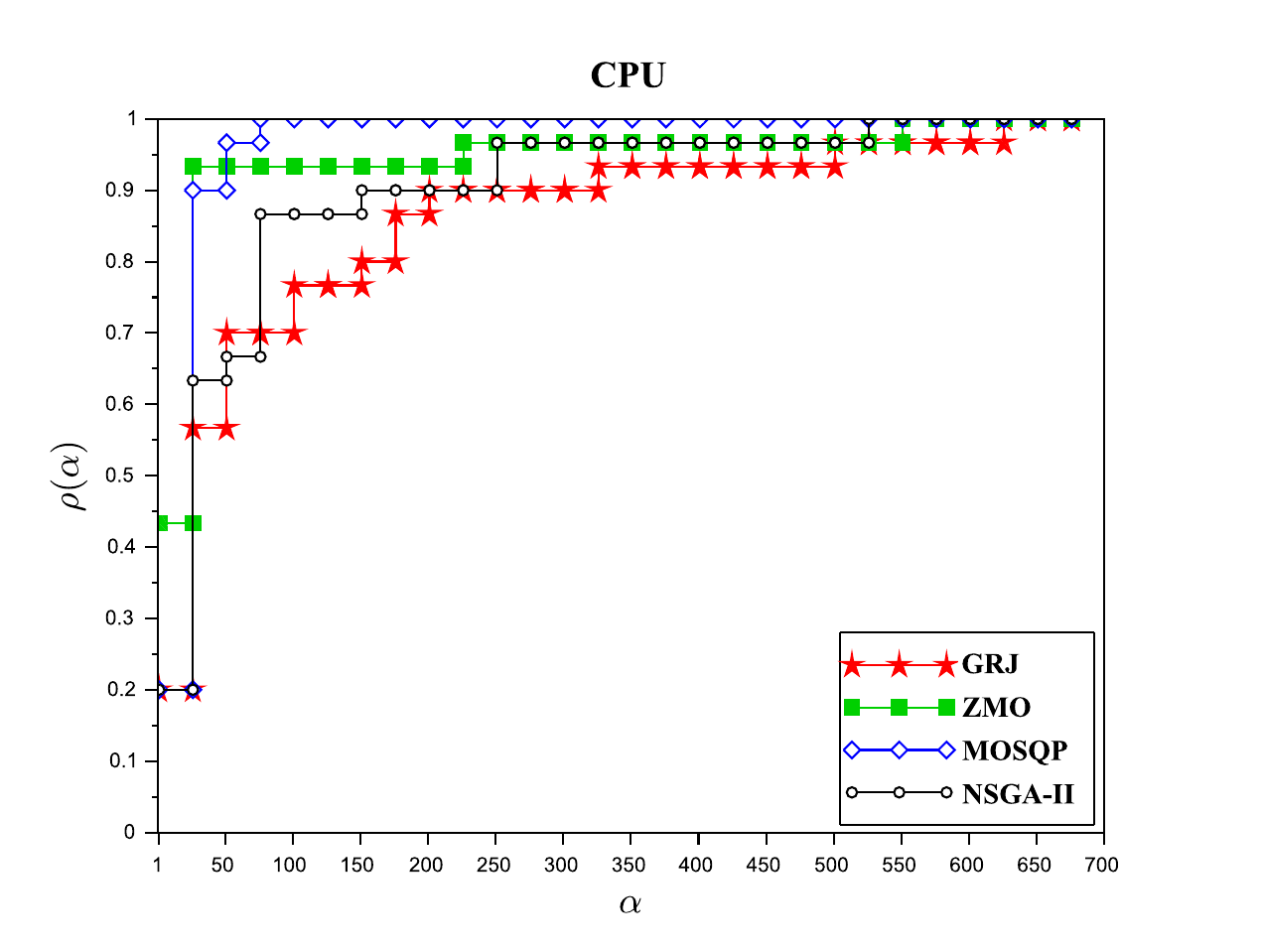}
	\end{subfigure}
	\hspace{0.6cm}
	\begin{subfigure}[b]{0.46\textwidth}
		\includegraphics[height=2.4in]{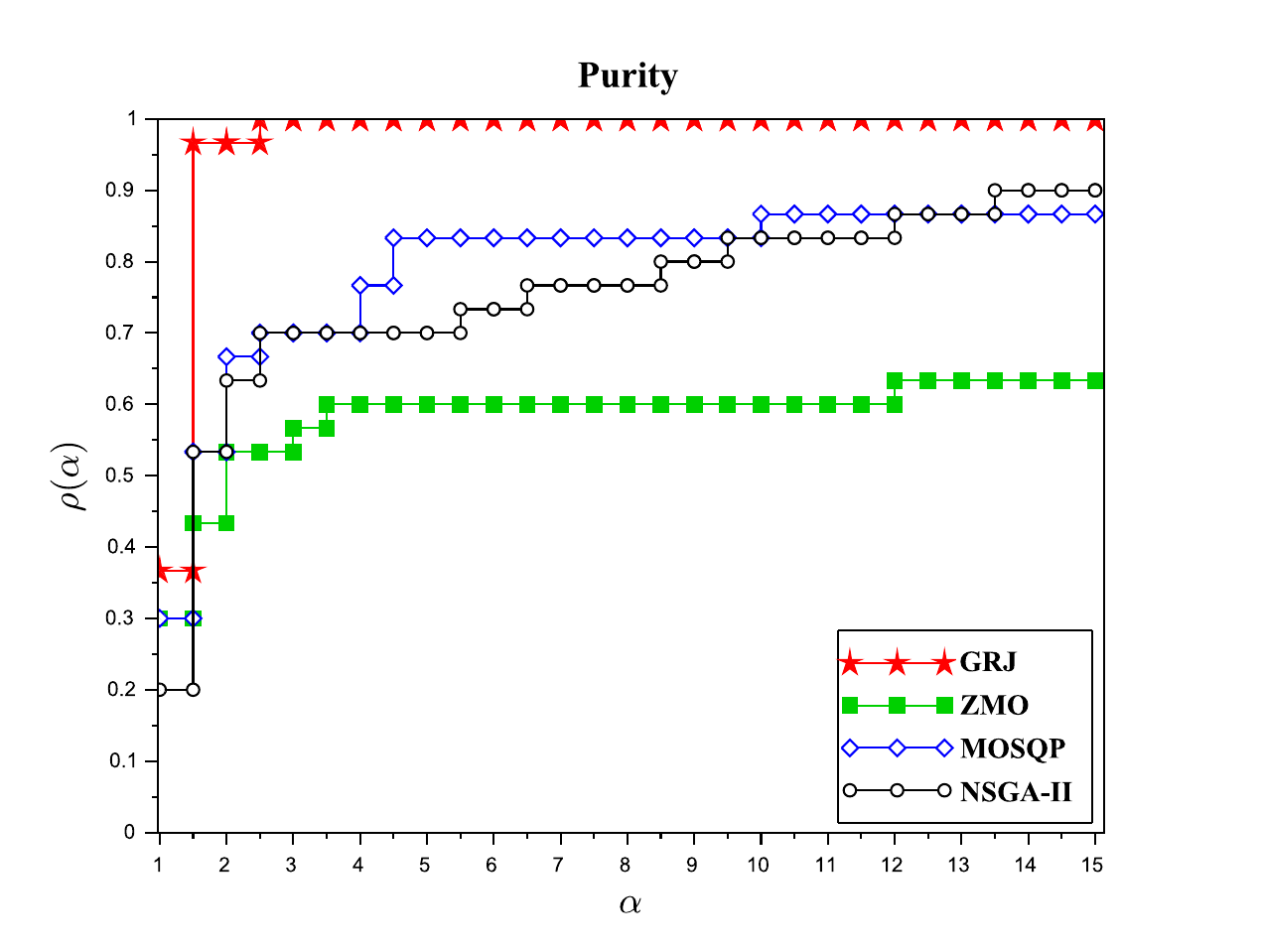}
	\end{subfigure}
	\vspace*{0.5cm}
	\hspace*{-0.8cm}	
	\begin{subfigure}[b]{0.46\textwidth}
		\includegraphics[height=2.4in]{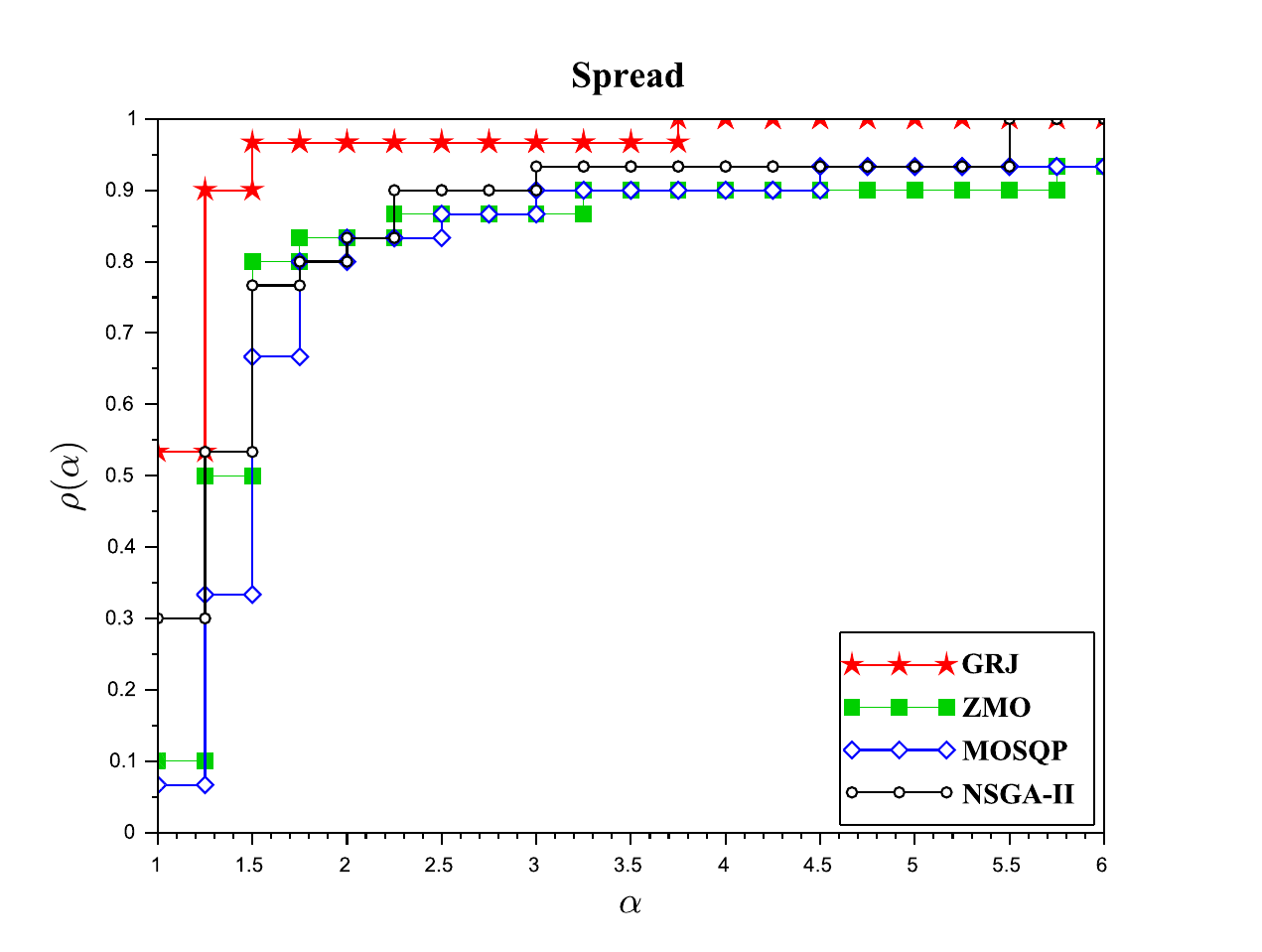}
	\end{subfigure}
	\hspace{0.6cm}
	\begin{subfigure}[b]{0.46\textwidth}
		\includegraphics[height=2.4in]{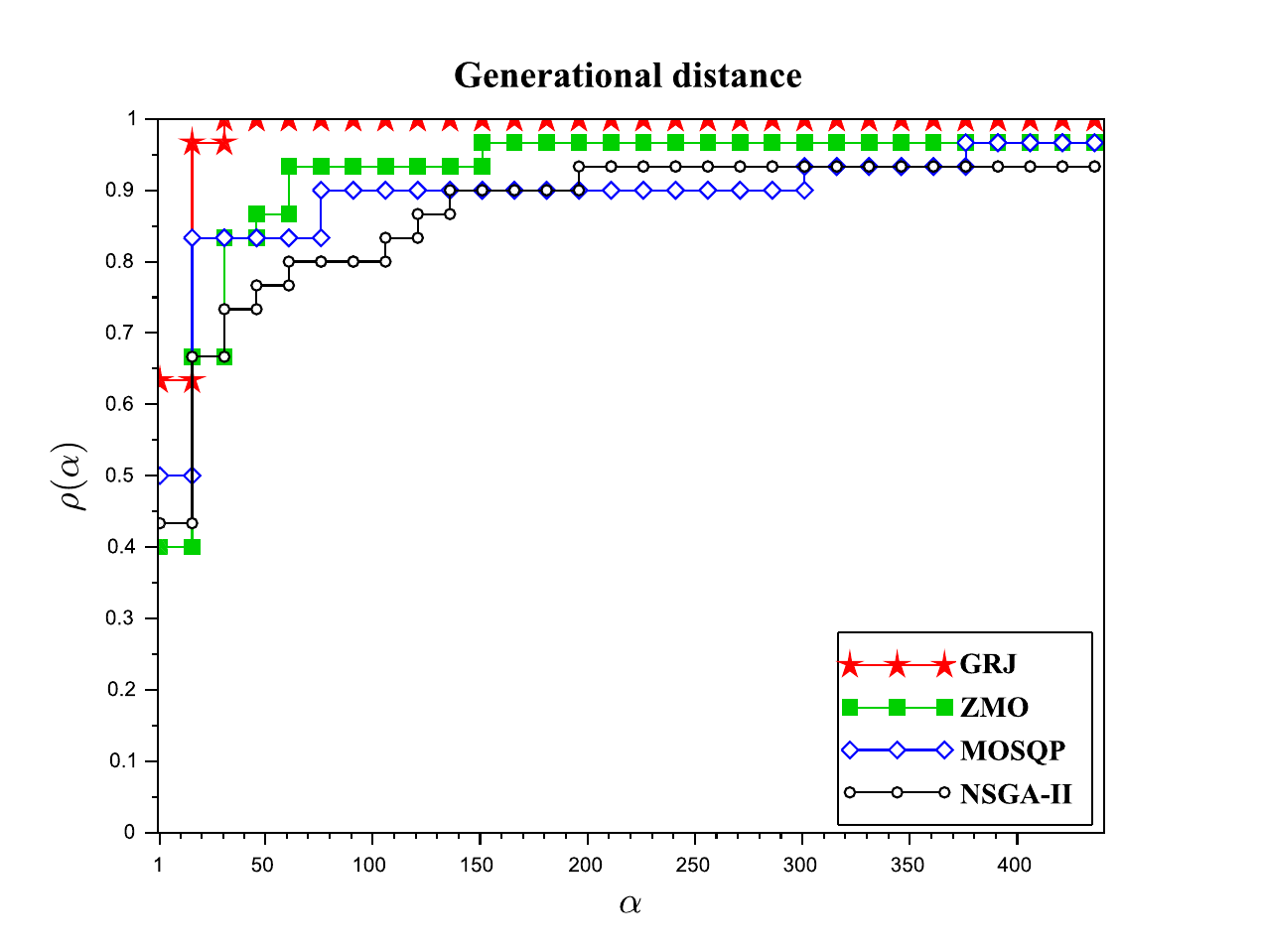}
	\end{subfigure}
    \vspace*{-0.5cm}
	\caption{Performance profiles}\label{fig1}
\end{figure}
  
\subsection{Real-world applications} 
As mentioned earlier, we conclude this section by providing additional details on the 
two real-world engineering problems: `Disc Brake' and `Welded Beam'. 

\subsubsection{`Disc Brake' design problem}

This problem consists to minimize simultaneously the mass of the brake and the stopping 
time. The variables represent the inner radius of the disc, the outer radius, the engaging 
force, and the number of friction surfaces. The constraints for the design include a minimum 
distance between the inner and outer radii, a maximum allowable brake length, and limitations 
related to pressure, temperature, and torque. The problem is a nonlinearly constrained bicriteria 
problem (see \cite{ray} for further details), formally defined as follows:
\begin{align*}\small
(\text{Disc Brake}) \qquad&\text{Minimize} \quad \left(4.9\times 10^{-5}\left(x_2^2-x_1^2\right)\left(x_4-1\right),\displaystyle{\frac{9.82\times 10^{6}\left(x_2^2-x_1^2\right)}{x_3x_4\left(x_2^3-x_1^3\right)}}\right),\\[0.2cm]
&\text{subject to}  \quad \left(x_2-x_1\right)+20\leq 0, \ \ 2.5\left(x_4+1\right)-30\leq 0,\\[0.2cm]
&\qquad \qquad \quad  \ \frac{x_3}{3.14\left(x_2^2-x_1^2\right)}-0.4\leq 0, \ \ \frac{2.22x_3\left(x_1^3-x_2^3\right)}{10^3\times\left(x_2^2-x_1^2\right)^2}-1\leq 0,\\[0.2cm] 
&\qquad \qquad \quad  \ \frac{2.66x_3x_4\left(x_1^3-x_2^3\right)}{10^2\left(x_2^2-x_1^2\right)^2}+900\leq 0,\\[0.2cm] 
&\qquad \qquad \quad  \ 55\leq x_1 \leq 80, \ \ 75\leq x_2 \leq 110,\\[0.2cm]
&\qquad \qquad \quad  \ 10^3\leq x_3 \leq 3\times 10^3, \ \ 2\leq x_4 \leq 20.
\end{align*}

\vspace*{-0.6cm}

\begin{figure}[H]
	\centering
	\includegraphics[height=4.2in]{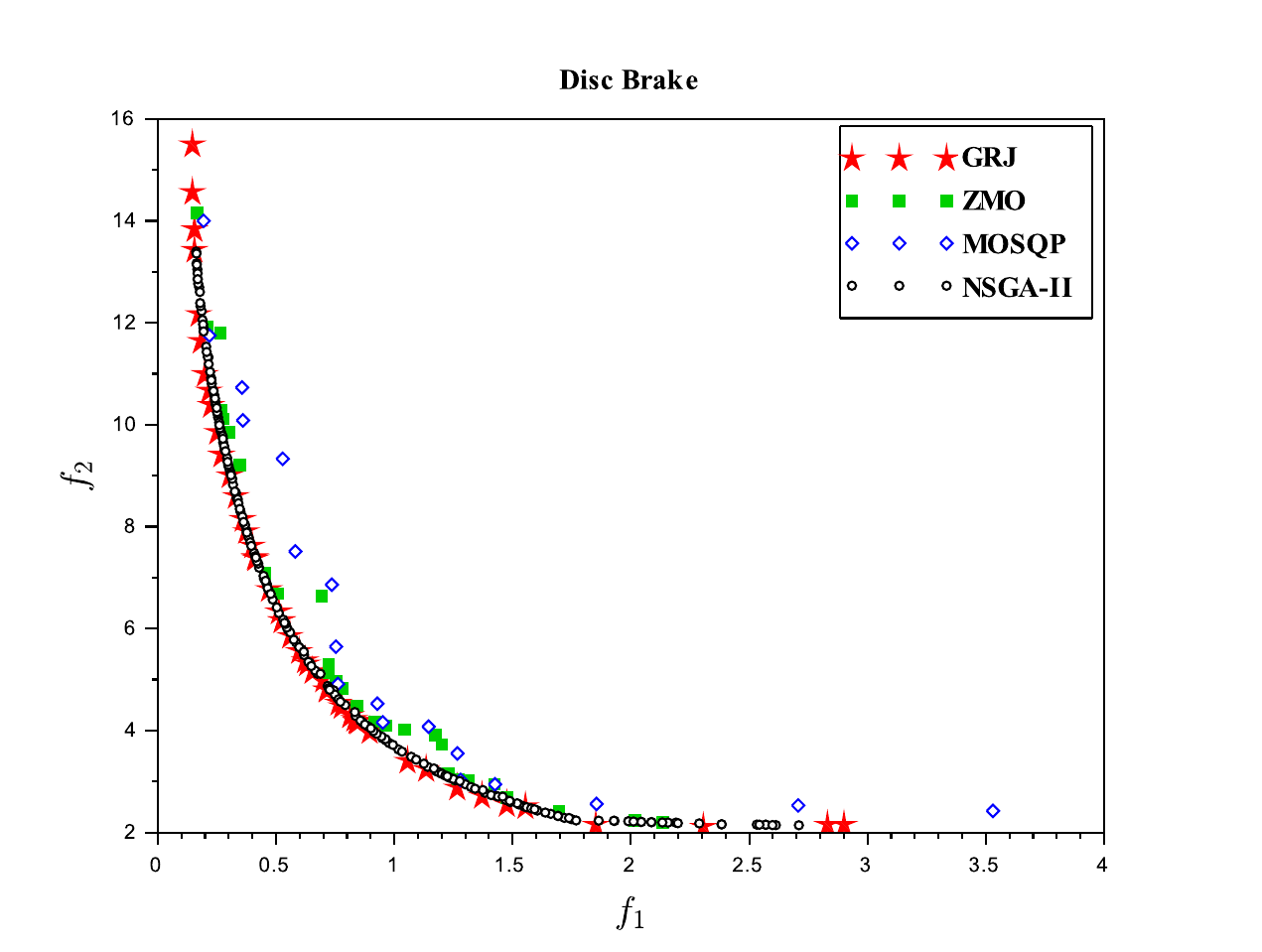}\vspace*{-0.2cm}
	\caption{Best Pareto front approximations for the `Disc Brake' design problem \\ 
		     by GRJ, ZMO, MOSQP and NSGA-II} \label{fig2}
\end{figure}

The graphical representation of the approximate Pareto fronts in Fig.~\ref{fig2} 
obtained by the four methods highlights the superior ability of GRJ and NSGA-II 
to better explore the Pareto front for the `Disc Brake' problem, in contrast to 
ZMO and MOSQP, which face serious challenges. Nevertheless, it is important to 
note that MOSQP slightly outperforms in this example ZMO, as it manages to find 
some relatively interesting solutions. 

\subsubsection{`Welded Beam' design problem}

This problem involves minimizing both the cost of fabrication and the deflection 
of the beam under an applied load subject to some constraints. The two objectives 
are inherently conflicting, since reducing deflection generally leads to higher 
manufacturing costs, primarily due to setup, material, and welding labor costs. 
The design involves four decision variables and four nonlinear constraints: shear 
stress, normal stress, weld length, and buckling limitations (see \cite{ray} for 
more details). Its formal definition is given as follows:
\begin{align*}\small
(\text{Welded Beam}) \qquad&\text{Minimize} \quad \left(1.10471x_1^2x_2+0.04811x_3x_4(14+x_3),\frac{2.1952}{x_3^3x_4}\right),\\[0.2cm]
&\text{subject to}  \quad \sqrt{\tau'(x)^2+\tau''(x)^2+\frac{x_2\tau'(x)\tau''(x)}{\sqrt{0.25\big(x_2^2+(x_1+x_3)^2\big)}}}\leq 13600,\\[0.2cm]
&\qquad \qquad \quad  \ \frac{1}{x_3^3x_4}\leq \frac{3}{504}, \ x_3x_4^3(0.0282346x_3-1)\leq \frac{6\times 10^3}{64746.022},\\[0.2cm]
&\qquad \qquad \quad  \ x_1\leq x_4, \ 0.125\leq x_1, x_4 \leq 5, \ \ 0.1\leq x_2, x_3 \leq 10,
\end{align*}

\medskip 
\noindent 
where $\tau'(x)=\frac{6\times 10^3}{\sqrt{2}x_1x_2}$ and $\tau''(x)=\frac{3\times 10^3\Big(14+0.5x_2\sqrt{0.25(x_2^2+(x_1+x_3)^2)}\Big)}{0.707x_1x_2\Big(\frac{x_2^2}{12}+0.25(x_1+x_3)^2\Big)}$.

\vspace*{-0.4cm}

\begin{figure}[H]
	\centering
	\includegraphics[height=4.2in]{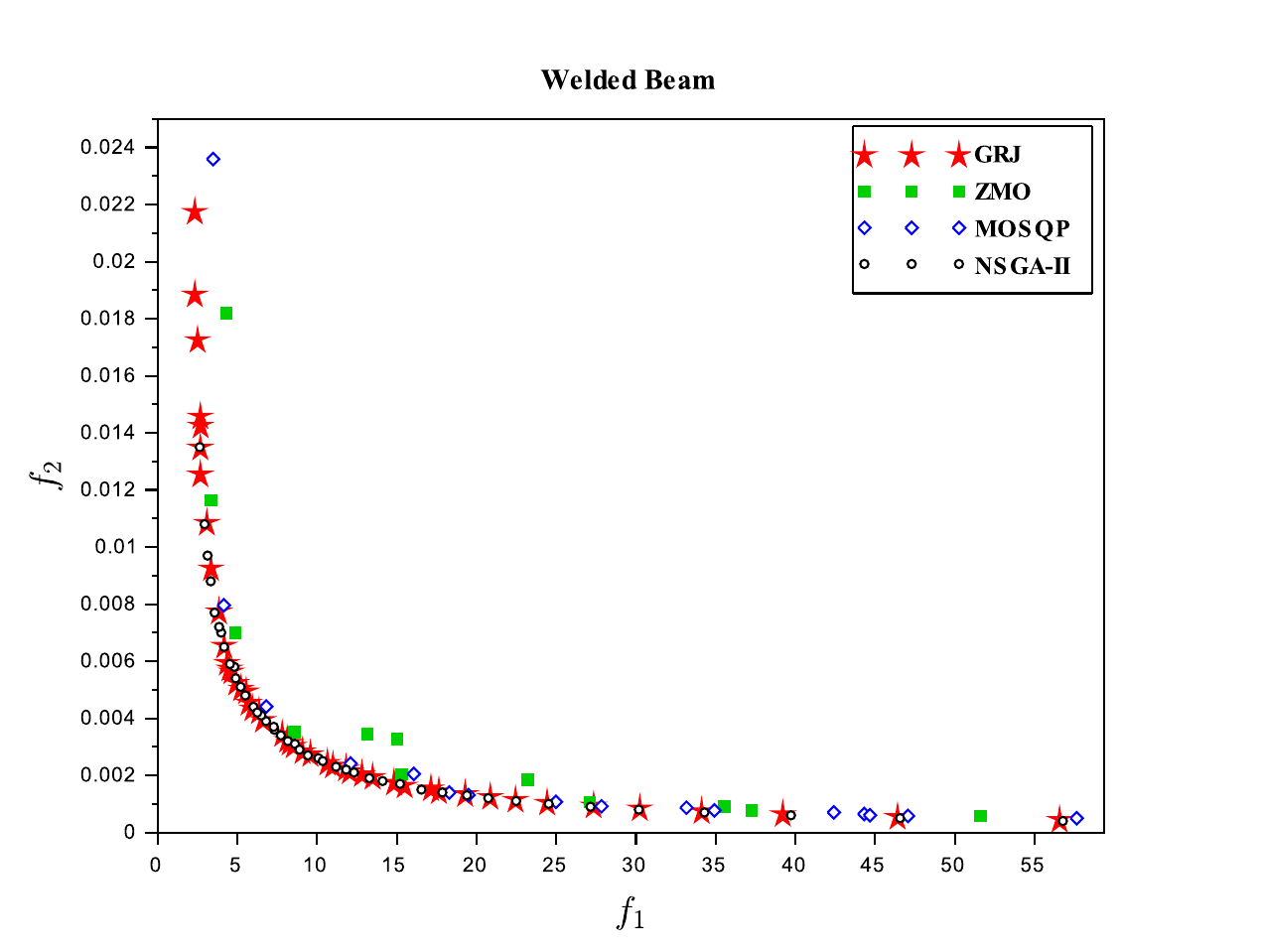}
	\vspace*{-0.2cm}
	\caption{Best Pareto front approximations for the `Welded Beam' design problem \\
		     by GRJ, ZMO, MOSQP and NSGA-II}\label{fig3}
\end{figure} 

In the context of the `Welded Beam' problem, as shown in Fig.~\ref{fig3}, 
all methods succeed in approximating significant regions of the Pareto front. 
However, the resulting dispersions differ from one method to another, so that 
none can be objectively favoured in this regard. In terms of convergence, it 
is also observed that the ZMO method is dominated by the other methods, which 
remain in strong competition with each other. This clearly demonstrates the 
effectiveness of our GRJ method in solving real-world problems. Moreover, the 
consistent performance of GRJ across various metrics highlights its robustness 
and reliability as a practical optimization approach.

\bigskip 
\noindent
{\footnotesize {\bf Statements and Declarations.} 
	The authors have no pertinent declarations concerning conflicts of interest, financial or 
	non-financial interests, competing interests or other statements to disclose.}

\end{document}